\documentclass[preprint]{elsarticle}

\usepackage{amsfonts}
\usepackage{booktabs}
\usepackage{color,colortbl} 
\usepackage{mathtools}
\usepackage{hyperref} 
\usepackage{multirow}
\usepackage{booktabs}
\usepackage{rotating}
\usepackage{graphicx,import} 
\usepackage{tikz}  
\usetikzlibrary{decorations.text}
\usepackage{caption} 

\DeclarePairedDelimiter{\ceil}{\lceil}{\rceil} 

\biboptions{authoryear}

\DeclareMathOperator{\cpoPresenceOf}        {presenceOf}
\DeclareMathOperator{\cpoEndOf}             {endOf}
\DeclareMathOperator{\cpoStartOf}           {startOf}
\DeclareMathOperator{\cpoSizeOf}            {sizeOf}
\DeclareMathOperator{\cpoNoOverlap}         {noOverlap}
\DeclareMathOperator{\cpoAlternative}       {alternative}
\DeclareMathOperator{\cpoSpan}              {span}
\DeclareMathOperator{\cpoEndBeforeStart}    {endBeforeStart}
\DeclareMathOperator{\cpoEndBeforeEnd}      {endBeforeEnd}
\DeclareMathOperator{\cpoEndAtStart}        {endAtStart}
\DeclareMathOperator{\cpoStartAtStart}      {startAtStart}
\DeclareMathOperator{\cpoInterval}          {interval}
\DeclareMathOperator{\cpoSequence}          {sequence}
\DeclareMathOperator{\cpoIntensity}         {intensity}
\DeclareMathOperator{\cpoSize}              {size}
\DeclareMathOperator{\cpoOptional}          {opt}
\DeclareMathOperator{\cpoPulse}             {pulse}

\DeclareMathOperator{\cpoSeqOn}             {on}
\DeclareMathOperator{\cpoSeqTypes}          {types}
\DeclareMathOperator{\cpoTypeOfPrev}        {typeOfPrev}
\DeclareMathOperator{\cpoForbidStart}       {forbidStart}
\DeclareMathOperator{\cpoForbidEnd}         {forbidEnd}
\DeclareMathOperator{\cpoForbidExtent}      {forbidExtent}

\begin{document}

\begin{frontmatter}

\title{Mixed Integer Linear Programming and Constraint Programming Models for the Online Printing Shop Scheduling Problem}

\author[unilu]{{Willian T. Lunardi}\corref{mycorrespondingauthor}}
\cortext[mycorrespondingauthor]{Corresponding author}
\author[imeusp]{{Ernesto G. Birgin}}
\author[ibm]{{Philippe Laborie}}
\author[poliusp]{{D\'ebora P. Ronconi}}
\author[unilu]{{Holger Voos}}

\address[unilu]{University of Luxembourg, 29 John F Kennedy, L-1855, Luxembourg, Luxembourg. e-mail: \{willian.tessarolunardi$\;|\;$holger.voos\}@uni.lu} 
\address[imeusp]{Department of Computer Science, Institute of Mathematics and Statistics, University of S\~ao Paulo, Rua do Mat\~ao, 1010, Cidade Universit\'aria, 05508-090, S\~ao Paulo, SP, Brazil. e-mail: egbirgin@ime.usp.br}
\address[ibm]{IBM France, 9 rue de Verdun, BP 85, 94253 Gentilly, France. e-mail: laborie@fr.ibm.com} 
\address[poliusp]{Department of Production Engineering, Polytechnic School, University of S\~ao Paulo, Av. Prof. Luciano Gualberto, 1380, Cidade Universit\'aria, 05508-010, S\~ao Paulo, SP, Brazil. e-mail: dronconi@usp.br}

\begin{abstract}
In this work, the online printing shop scheduling problem is considered. This challenging real problem, that appears in the nowadays printing industry, can be seen as a flexible job shop scheduling problem with sequence flexibility in which precedence constraints among operations of a job are given by an arbitrary directed acyclic graph. In addition, several complicating particularities such as periods of unavailability of the machines, resumable operations, sequence-dependent setup times, partial overlapping among operations with precedence constraints, release times, and fixed operations are also present in the addressed problem. In the present work, mixed integer linear programming and constraint programming models for the minimization of the makespan are presented. Modeling the problem is twofold. On the one hand, the problem is precisely defined. On the other hand, the capabilities and limitations of a commercial software for solving the models are analyzed. Extensive numerical experiments with small-, medium-, and large-sized instances are presented. Numerical experiments show that the commercial solver is able to optimally solve only a fraction of the small-sized instances when considering the mixed integer linear programming model; while all small-sized and a fraction of the medium-sized instances are optimally solved when considering the constraint programming formulation of the problem. Moreover, the commercial solver is able to deliver feasible solutions for the large-sized instances that are of the size of the instances that appear in practice.
\end{abstract}

\begin{keyword}
Flexible job shop scheduling with sequence flexibility \sep Resumable operations \sep Unavailability of the machines \sep Sequence-dependent setup time \sep Mixed integer linear programming \sep Constraint programming\\[2mm]

\MSC[2010] 90B35 \sep 90C11 \sep  90C59
\end{keyword}

\end{frontmatter}

\section{Introduction}\label{sec:intro}

The online printing shop scheduling problem emerged in the printing industry that nowadays faces challenges that were unheard of just a couple of years ago. Printing jobs are getting more complex; while demand is decreasing due to the expanding use of online advertising and the decreasing of the printed matter. Online printing shops (OPS) commenced taking advantage of the Internet and standardization towards mass customization to attract more clients and to reduce costs. This new business model enables printers to leave the inquiry and order process to the customer. Alternatively to a small number of fixed contracts, online printing shops can reach up to 20k online orders of diverse clients per day. Through an online system, there are placed orders that can be of various types, such as bags, beer mats, books, brochures, business cards, calendars, certificates, Christmas cards, envelopes, flyers, folded leaflets, greeting cards, multimedia CD/DVD, napkins, paper cups, photo books, posters, ring binders, stamps, stickers, tags, tickets, wristbands, among others. The system guides the user through a standardized process that defines details including size, orientation, material, finish, color mode, design, among others. 

\begin{figure}[ht!]
\centering
\resizebox{0.9\textwidth}{!}{
\begin{tabular}{ccc}
\includegraphics{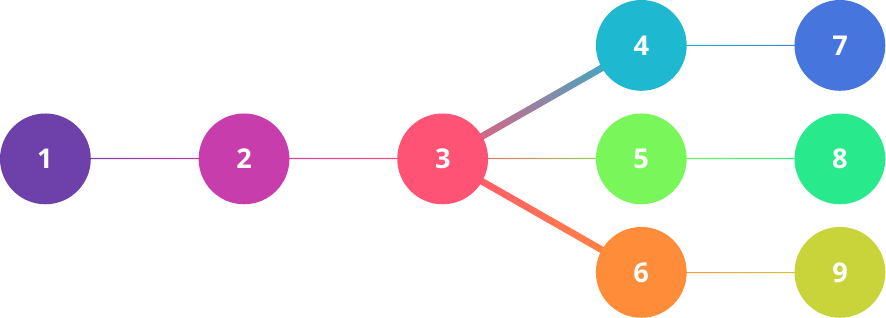} & \quad \quad \quad \quad &
\includegraphics{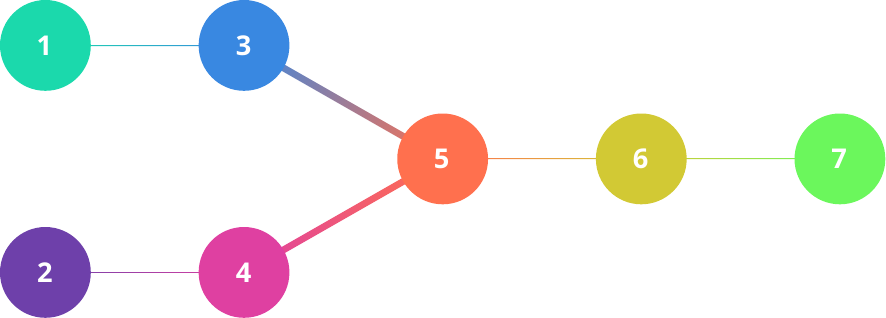} \\[4mm]
(a) Example of a job with a ganging operation. & \quad \quad \quad \quad &
(b) Example of a job with a binding operation.
\end{tabular}}
\caption{Directed acyclic graphs representing the precedence constraints between operations of two jobs. Each node represents an operation and arcs are directed from left to right.}
\label{fig:ops_job}
\end{figure}

\begin{figure}[ht!]
\centering
\footnotesize
\includegraphics[width=\textwidth]{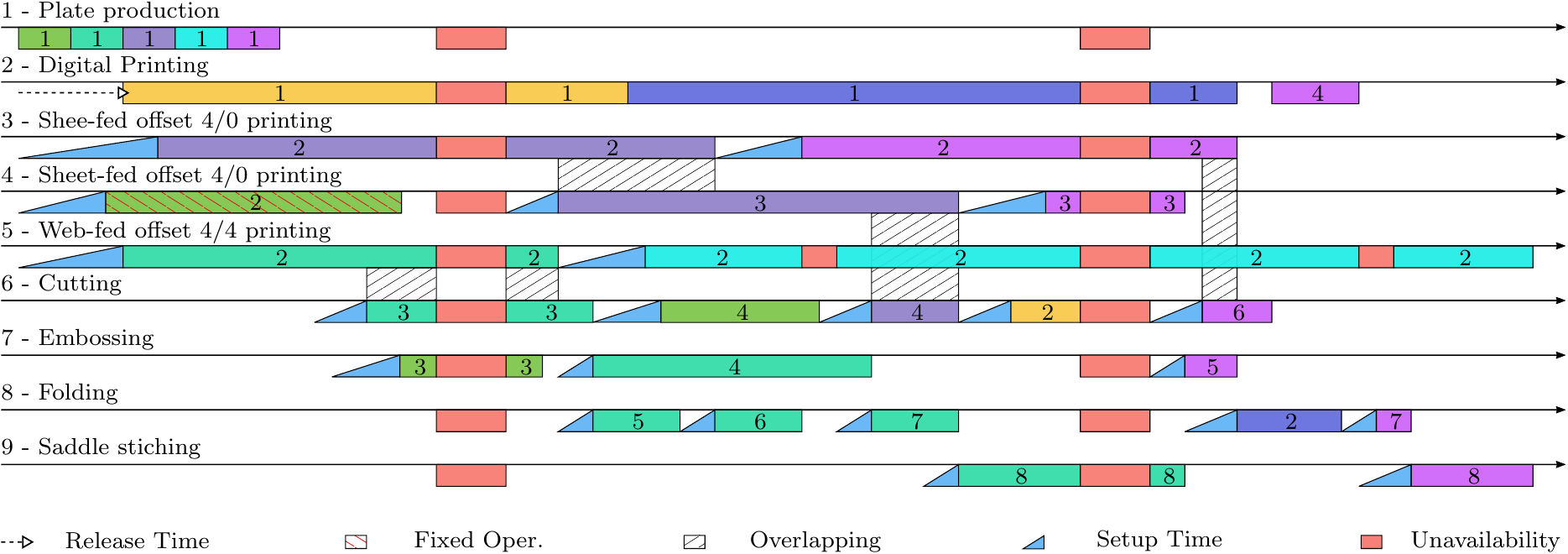}
\caption{Gantt chart of a solution to a simple instance of the considered scheduling problem containing~9 machines and~7 jobs. The 4/0 in the machine description means it is capable of printing 4 colors one side; while 4/4 means 4 colors both sides. Digital printing operations do not require plate production. Disassembling and assembling operations may occur. The jobs' routes through the machines are not the same. Operations of a job are represented with the same color.}
\label{fig:ops_gantt}
\end{figure}

Clearly, products to be manufactured have a different production plan; and the production of all of them involves a printing operation. When a significant number of orders is reached, looking at the printing operation of each product and aiming to reduce the production cost, a cutting stock problem is solved to join the printing operations of different placed orders. These merged printing operations are known as \textit{ganging operations}. Orders whose printing operations are ganged constitute a single job. Thus, jobs of the considered scheduling problem, composed by a heterogeneous set of operations with arbitrary precedence constraints, have wide diversity. Operations of a job can be categorized into three major groups: \textit{prepress}, the preparation of materials and plates for printing; \textit{press}, the actual printing process; and \textit{postpress}, the cutting, folding, binding, embossing (varnishing, laminating, hot foil, etc.), and trimming of printed sheets into their final form. However, it is important to highlight that not all jobs have to follow the same route, e.g.,\ not all jobs need a cutting operation, not all jobs present folding operations, not all jobs present embossing operations, etc. Moreover, the order in which operations of a job take place differs from job to job, e.g.,\ a job may require a printing operation followed by a cutting operation and another printing operation. Another job may require embossing before cutting; while another job may require embossing after cutting or it may require several printing operations in sequence. Operations can be processed by several machines with varying processing times. Figure~\ref{fig:ops_job} schematically shows two possible job topologies. The figure aims to stress that, differently from the classical linear order, operations of a job have arbitrary precedence constraints represented by a directed acyclic graph (DAG). The DAG presented in Figure~\ref{fig:ops_job}a represents the topology of a job with a ganging operation. The ganging operation (node~2) is disassembled by a cutting operation (node~3) into three independent sequences. The topology presented in Figure~\ref{fig:ops_job}b represents, in a very simplified way, the production of a book. Independent sequences represent the parallel production of the book cover (nodes~1 and~3) and book pages (nodes~2 and~4), later assembled by a binding operation (node~5). Figure~\ref{fig:ops_gantt} presents a Gantt chart of a toy instance of the OPS scheduling problem. The figure aims to illustrate a small set of diverse jobs in a complicate system that includes features such as periods of unavailability of the machines, release times, setup operations, overlapping between operations of a job with precedence constraints, etc. 

In the OPS scheduling problem, a precedence relation saying that an operation~$i$ must be processed before an operation~$j$ is represented by an arc $(i,j)$ in a directed acyclic graph. On the one hand, operations with a precedence relation may overlap. The extent of this overlapping is limited by a fraction that defines the proportion of operation~$i$ that has to be processed before its successor~$j$ can start. On the other hand, the precedence implies that~$j$ can not be completed before~$i$. Several flexible machines are available in the OPS. Each operation~$i$ can be processed by a subset of machines, with varying processing times. Before starting processing, the machine must be configured according to the characteristics of each operation such as color, size, grammage, scheme, etc. The setup time is the time taken to set up a machine from the current setting to the setting of the next operation to be processed. The length of the setup time depends on the similarities between the two consecutive operations. The more different the operations' settings, the longer the machine setup. Another important feature of the considered problem is that machines may be unavailable during some planned (i.e.,\ known beforehand) time periods, being a consequence of pre-scheduling, maintenance, shift patterns, or overlapping of two consecutive time planning horizons. Operations can be interrupted by one or more unavailable periods, i.e.,\ operations are resumable. On the other hand, the machine setup cannot be interrupted by unavailable periods; and the end of the setup must be immediately followed by the beginning of the operation processing. A note about the interruption of the operations' processing by periods of unavailability of the machines is in order. According to Pinedo~\cite[p.16]{pinedo}, preemption implies that it is not necessary to keep an operation on a machine, once started, until its completion. The scheduler is allowed to interrupt the processing of an operation (preempt) at any point in time and put a different operation on the machine instead. This is not exactly the case of the OPS scheduling problem. Preemption of this kind is not allowed. In the OPS scheduling problem, when a machine starts the execution of an operation, the operation remains on the machine until its completion. This means that the processing may be interrupted by a period of unavailability of the machine, but the processing of the operation must be resumed as soon as the machine returns to be available, i.e.,\ the operation must remain on the machine until completion. Also in order is a note about the fact that setup times can not be interrupted and that they must be immediately followed by the execution of an operation, with no idle time between them. A setup operation could be the operation of cleaning a machine; while a period of unavailability of a machine could correspond to pre-scheduled maintenance. The machine can not be opened and half-cleaned, the maintenance operation executed, and the cleaning operation completed. This is because the maintenance operation may ruin the cleaning operation. The same situation occurs if the period of unavailability corresponds to a night shift during which the store is closed. In this case, the half-cleaned opened machine could get dirty because of dust or insects during the night. A customer may request to visit the shop to check her/his order while it is being produced. To meet with the customer's visiting schedule, some operations may be fixed at particular resources with a pre-defined starting time. During the scheduling optimization, fixed operations cannot be moved or reassigned. Fixed operations are also useful to include, in the scheduling to be constructed, operations of previous scheduling that are still being processed in the system. Considering the (final) operations of the previous scheduling as fixed operations of the new one allows one to adequately take into account the corresponding sequence-dependent setup times.

Based on the manufacturing information system and key features of the production line, a production schedule should be planned to maximize the production effectiveness so that online printing shops can gain as much production benefit as possible. The production effectiveness can be represented by an objective function under the constraint of manufacturing resources (machines, information sources, production workflows, etc.) A usual objective is to minimize the time that elapses from the start of work to the end, as in this way enterprises can reduce the cost of labor and electricity while achieving a high quality of products. In this paper, we study the problem of minimizing the makespan of an OPS scheduling problem by considering it as a flexible job shop (FJS) scheduling problem in which operations' precedence constraints are given by a DAG. The flexibility of representing the precedences with a DAG instead of with a linear order (i.e. a single sequence) fits within the concept of sequence flexibility. Note that representing the precedences with a DAG also allows some operations of a job to be processed in parallel. Resumable operations, periods of unavailability of the machines, sequence-dependent setup times, partial overlapping among operations with precedence constraints, release times, and fixed operations are additional features of the considered scheduling problem.

\section{Literature Review}
In the classical FJS scheduling problem, each operation can be processed by a subset of machines, situation known as routing flexibility. A wide range of techniques have been presented in the literature to deal with the FJS scheduling problem. A mixed integer linear programming (MILP) model for this problem is given in~\cite{fattahi2007mathematical}. An alternative model for the FJS scheduling problem is given in~\cite{ozguven2010mathematical}. This model also considers the flexibility of the job process plans -- while in the FJS problem a job is given by a set of operations with linear-order precedence constraints; the process plan flexibility means that each job can be completed by performing one among several alternative sets of operations, each one with linear-order precedence constraints. A comparison among MILP models for the FJS scheduling problem is presented in~\cite{demir2013evaluation}. Constraint programming (CP) formulations for the FJS scheduling problem are presented in~\cite{vilim2015failure} and~\cite{Laborie2018a}. Numerical results in these works show that the IBM ILOG CP Optimizer~\citep{Laborie2018} improves best-known results for several classical FJS instances, suggesting that this off-the-shelf solver can compete with or even outperform specialized algorithms. An extensive literature review of formulations and exact and heuristic methods developed in the last three decades to approach the FJS scheduling problem and some extensions is presented in~\cite{chaudhry2016research}. On the other hand, only a few works, mostly inspired by practical applications from the glass, mold, and printing industries, deal with the FJS scheduling problem with sequence flexibility. The literature review below, presented in chronological order, focuses on models and practical applications of the FJS scheduling problem with sequence flexibility. It aims to show that no paper in the literature proposes a model for or describes a practical application of the FJS scheduling problem with sequence flexibility encompassing all the challenging and complicating features of the OPS scheduling problem. Moreover, up to the authors' acknowledge, no paper in the literature tackles this problem with constraint programming techniques.

In~\cite{gan}, an FJS scheduling problem with sequence flexibility and process plan flexibility is considered. The problem is based on a practical application in a mold manufacturing shop; and it is tackled with a branch and bound algorithm. The simultaneous optimization of the process plan and the scheduling problem is uncommon in the literature. On the other hand, no model is presented and none of the additional features of the OPS scheduling problem appear in the mold manufacturing shop application. In~\cite{kim}, the integrated problem considered in~\cite{gan} is also addressed; and a symbiotic evolutionary algorithm is proposed. In comparison to the problem addressed in the present work, the problem appears not to be directly related to a practical application, no model is proposed, and none of the additional features of the OPS scheduling problem are considered. In~\cite{alvarez2005heuristic}, an FJS scheduling problem with sequence flexibility issued from the glass industry is addressed. A heuristic algorithm combining priority rules and local search is proposed. The considered problem includes extensions such as overlapping among operations, operations that can be processed simultaneously, and fixed intervals of unavailability of the machines. However, the problem does not include all the OPS scheduling problem's features and no model is given. In~\cite{vilcot2008tabu}, a problem issued from the printing industry, where some operations can be carried out at the same time, is addressed. In the considered problem, operations’ precedence constraints are such that each operation may have one or more predecessor but it can not have more than one successor, i.e.,\ ganging operations like the ones depicted in Figure~\ref{fig:ops_job}a are not considered. A bi-objective genetic algorithm based on the NSGA II is proposed to tackle the problem. The problem comes from the printing industry, as the OPS scheduling problem, and, in consequence, it includes some of the features of the OPS scheduling problem. However, no model is presented and a limited type of sequence flexibility is considered. In~\cite{sanghyup}, an FJS scheduling problem with AND/OR precedence constraints in the operations is considered. Thus, the problem corresponds to an FJS scheduling problem with sequence flexibility and process plan flexibility as the one considered in~\cite{kim} and in~\cite{gan}. A MILP model and genetic and tabu search algorithms are proposed. Release times for the jobs are considered, but none of the other additional features of the OPS scheduling problem are addressed. In~\cite{birgin2014milp}, a MILP model for the FJS scheduling problem in which precedence constraints between operations are described by an arbitrary directed acyclic graph is introduced; and the model for the FJS scheduling problem introduced in~\cite{ozguven2010mathematical} is extended to include sequence flexibility. List scheduling and beam search algorithms for the same problem are introduced in~\cite{birgin2015list}. None of the additional features of the OPS scheduling problem are considered in~\cite{birgin2014milp,birgin2015list}, where the considered problem is a simplification of a real problem coming from the printing industry. The scheduling of repairing orders and allocation of workers in an automobile collision repair shop is addressed in~\cite{andradescheduling}. The underlying scheduling problem is a dual-resource flexible job shop scheduling problem with sequence flexibility that aims to minimize a combination of makespan and mean tardiness. The authors extend the MILP formulation proposed in~\cite{birgin2014milp} and propose a constructive iterated greedy heuristic. Table~\ref{tabliterature} summarizes the main features of the FJS scheduling problem with sequence flexibility addressed in each one of the cited references.

\begin{table}[ht!]
\begin{center}
\begin{tabular}{|c|ccc|ccc|cccccc|}
\hline
 & \rotatebox{90}{Practical application} 
 & \rotatebox{90}{MILP formulation} 
 & \rotatebox{90}{CP formulation} 
 & \rotatebox{90}{Routing flexibility} 
 & \rotatebox{90}{Sequence flexibility} 
 & \rotatebox{90}{Process plan flexibility} 
 & \rotatebox{90}{Sequence-dependent setups $\phantom{a}$} 
 & \rotatebox{90}{Partial overlapping} 
 & \rotatebox{90}{Resumable operations} 
 & \rotatebox{90}{Machines' unavailabilities}  
 & \rotatebox{90}{Release times}  
 & \rotatebox{90}{Fixed operations} \\
\hline
\hline
\cite{gan} &
\multicolumn{1}{c!{\color{lightgray}\vrule}}{$\checkmark$} &
\multicolumn{1}{c!{\color{lightgray}\vrule}}{} &
&
\multicolumn{1}{c!{\color{lightgray}\vrule}}{$\checkmark$} &
\multicolumn{1}{c!{\color{lightgray}\vrule}}{$\checkmark$} &
$\checkmark$ &
\multicolumn{1}{c!{\color{lightgray}\vrule}}{} &
\multicolumn{1}{c!{\color{lightgray}\vrule}}{} &
\multicolumn{1}{c!{\color{lightgray}\vrule}}{} &
\multicolumn{1}{c!{\color{lightgray}\vrule}}{} &
\multicolumn{1}{c!{\color{lightgray}\vrule}}{} & \\
\arrayrulecolor{lightgray}
\cline{2-13}
\cite{kim} &
\multicolumn{1}{c!{\color{lightgray}\vrule}}{} &
\multicolumn{1}{c!{\color{lightgray}\vrule}}{} &
&
\multicolumn{1}{c!{\color{lightgray}\vrule}}{$\checkmark$} &
\multicolumn{1}{c!{\color{lightgray}\vrule}}{$\checkmark$} &
$\checkmark$ &
\multicolumn{1}{c!{\color{lightgray}\vrule}}{} &
\multicolumn{1}{c!{\color{lightgray}\vrule}}{} &
\multicolumn{1}{c!{\color{lightgray}\vrule}}{} &
\multicolumn{1}{c!{\color{lightgray}\vrule}}{} &
\multicolumn{1}{c!{\color{lightgray}\vrule}}{} & \\
\cline{2-13}
\cite{alvarez2005heuristic} &
\multicolumn{1}{c!{\color{lightgray}\vrule}}{$\checkmark$} &
\multicolumn{1}{c!{\color{lightgray}\vrule}}{} &
&
\multicolumn{1}{c!{\color{lightgray}\vrule}}{$\checkmark$} &
\multicolumn{1}{c!{\color{lightgray}\vrule}}{$\checkmark$} &
&
\multicolumn{1}{c!{\color{lightgray}\vrule}}{} &
\multicolumn{1}{c!{\color{lightgray}\vrule}}{$\checkmark$} &
\multicolumn{1}{c!{\color{lightgray}\vrule}}{$\checkmark$} &
\multicolumn{1}{c!{\color{lightgray}\vrule}}{$\checkmark$} &
\multicolumn{1}{c!{\color{lightgray}\vrule}}{} & \\
\cline{2-13}
\cite{vilcot2008tabu} &
\multicolumn{1}{c!{\color{lightgray}\vrule}}{$\checkmark$} &
\multicolumn{1}{c!{\color{lightgray}\vrule}}{} &
&
\multicolumn{1}{c!{\color{lightgray}\vrule}}{$\checkmark$} &
\multicolumn{1}{c!{\color{lightgray}\vrule}}{$\checkmark$} &
&
\multicolumn{1}{c!{\color{lightgray}\vrule}}{} &
\multicolumn{1}{c!{\color{lightgray}\vrule}}{} &
\multicolumn{1}{c!{\color{lightgray}\vrule}}{} &
\multicolumn{1}{c!{\color{lightgray}\vrule}}{} &
\multicolumn{1}{c!{\color{lightgray}\vrule}}{} & \\
\cline{2-13}
\cite{sanghyup} &
\multicolumn{1}{c!{\color{lightgray}\vrule}}{$\checkmark$} &
\multicolumn{1}{c!{\color{lightgray}\vrule}}{$\checkmark$} &
&
\multicolumn{1}{c!{\color{lightgray}\vrule}}{$\checkmark$} &
\multicolumn{1}{c!{\color{lightgray}\vrule}}{$\checkmark$} &
$\checkmark$ &
\multicolumn{1}{c!{\color{lightgray}\vrule}}{} &
\multicolumn{1}{c!{\color{lightgray}\vrule}}{} &
\multicolumn{1}{c!{\color{lightgray}\vrule}}{} &
\multicolumn{1}{c!{\color{lightgray}\vrule}}{} &
\multicolumn{1}{c!{\color{lightgray}\vrule}}{$\checkmark$} & \\
\cline{2-13}
\cite{birgin2014milp} &
\multicolumn{1}{c!{\color{lightgray}\vrule}}{} &
\multicolumn{1}{c!{\color{lightgray}\vrule}}{$\checkmark$} &
&
\multicolumn{1}{c!{\color{lightgray}\vrule}}{$\checkmark$} &
\multicolumn{1}{c!{\color{lightgray}\vrule}}{$\checkmark$} &
&
\multicolumn{1}{c!{\color{lightgray}\vrule}}{} &
\multicolumn{1}{c!{\color{lightgray}\vrule}}{} &
\multicolumn{1}{c!{\color{lightgray}\vrule}}{} &
\multicolumn{1}{c!{\color{lightgray}\vrule}}{} &
\multicolumn{1}{c!{\color{lightgray}\vrule}}{} & \\
\cline{2-13}
\cite{birgin2015list} &
\multicolumn{1}{c!{\color{lightgray}\vrule}}{} &
\multicolumn{1}{c!{\color{lightgray}\vrule}}{$\checkmark$} &
&
\multicolumn{1}{c!{\color{lightgray}\vrule}}{$\checkmark$} &
\multicolumn{1}{c!{\color{lightgray}\vrule}}{$\checkmark$} &
&
\multicolumn{1}{c!{\color{lightgray}\vrule}}{} &
\multicolumn{1}{c!{\color{lightgray}\vrule}}{} &
\multicolumn{1}{c!{\color{lightgray}\vrule}}{} &
\multicolumn{1}{c!{\color{lightgray}\vrule}}{} &
\multicolumn{1}{c!{\color{lightgray}\vrule}}{} & \\
\cline{2-13}
\cite{andradescheduling} &
\multicolumn{1}{c!{\color{lightgray}\vrule}}{$\checkmark$} &
\multicolumn{1}{c!{\color{lightgray}\vrule}}{$\checkmark$} &
&
\multicolumn{1}{c!{\color{lightgray}\vrule}}{$\checkmark$} &
\multicolumn{1}{c!{\color{lightgray}\vrule}}{$\checkmark$} &
&
\multicolumn{1}{c!{\color{lightgray}\vrule}}{} &
\multicolumn{1}{c!{\color{lightgray}\vrule}}{} &
\multicolumn{1}{c!{\color{lightgray}\vrule}}{} &
\multicolumn{1}{c!{\color{lightgray}\vrule}}{} &
\multicolumn{1}{c!{\color{lightgray}\vrule}}{} & \\
\cline{2-13}
Current work &
\multicolumn{1}{c!{\color{lightgray}\vrule}}{$\checkmark$} &
\multicolumn{1}{c!{\color{lightgray}\vrule}}{$\checkmark$} &
$\checkmark$ &
\multicolumn{1}{c!{\color{lightgray}\vrule}}{$\checkmark$} &
\multicolumn{1}{c!{\color{lightgray}\vrule}}{$\checkmark$} &
&
\multicolumn{1}{c!{\color{lightgray}\vrule}}{$\checkmark$} &
\multicolumn{1}{c!{\color{lightgray}\vrule}}{$\checkmark$} &
\multicolumn{1}{c!{\color{lightgray}\vrule}}{$\checkmark$} &
\multicolumn{1}{c!{\color{lightgray}\vrule}}{$\checkmark$} &
\multicolumn{1}{c!{\color{lightgray}\vrule}}{$\checkmark$} &
$\checkmark$ \\
\arrayrulecolor{black}
\hline
\end{tabular}
\end{center}
\caption{Summary of the literature review considering works that deal with formulations and/or practical applications of the FJS scheduling problem with sequence flexibility. (By definition, all mentioned works have a checkmark in the columns related to routing and sequence flexibilities.)}
\label{tabliterature}
\end{table}

The purpose of this work is to tackle a challenging real scheduling problem; so presenting the problem formulation in a precise and universal language (MILP formulation) is the first goal of the work. Since the problem is NP-hard -- it is a generalization of the job shop scheduling problem, known to be NP-hard~\citep{nphard}, it is well-known that only small-sized instances can be solved with a certificate of optimality. Numerical experiments in this work show that this is exactly the case when the MILP formulation is tackled with an exact commercial solver~\citep{cplexmanual}. On the other hand, the commercial solver IBM ILOG CP Optimizer~\citep{Laborie2018}, that applies to CP formulations, is presented in the literature~\citep{vilim2015failure} as a half-heuristic-half-exact method that, when it is not able to produce a solution with a guarantee of optimality, it produces a good quality solution in a reasonable amount of time. Verifying the size of the instances that IBM ILOG CP Optimizer can solve to optimality in comparison to IBM ILOG CPLEX is a second goal of the manuscript. As a side effect, both models (MILP and CP) can be validated and analyzed comparatively; and a benchmark set with known optimal solutions can be built. As a third goal, this works aims to analyze the capability of IBM ILOG CP Optimizer to find good quality solutions when applied to large-sized instances that are of the size of the instances of the OPS scheduling problem that appear in practice. As a whole, it is expected the present work to provide in-depth knowledge of the OPS scheduling problem that leverages the future development of ad-hoc heuristic methods to be applied in practice.

The paper is structured as follows. In Section~\ref{sec:mip}, the considered problem is described in detail; and its precise formulation is presented through a MILP formulation. Anticipating that, due to its complexity, only small-sized instances would be solved applying a commercial solver to the MILP formulation, a CP formulation is presented in Section~\ref{sec:cp}. A short discussion comparing both models is given in Section~\ref{sec:disc}. In Section~\ref{sec:exp}, extensive numerical experiments are presented. First, small-sized instances of the MILP and CP Optimizer formulations are solved using IBM ILOG CPLEX~\citep{cplexmanual} and IBM ILOG CP Optimizer~\citep{Laborie2018}, respectively. Then, in order to assert the capacity of the CP Optimizer solver for finding optimal solutions, medium-sized instances are considered; and several alternative models and resolution strategies are analyzed. Finally, in order to evaluate the applicability of the CP Optimizer solver in reality, numerical experiments are conducted with large-sized instances. Conclusions and lines for future research are given in Section~\ref{sec:concl}.

\section{Mixed integer linear programming formulation}\label{sec:mip}

The FJS scheduling problem with sequence flexibility considered in~\cite{birgin2014milp} generalizes the classical FJS scheduling problem in the sense that the linear order of the operations of a job is replaced by arbitrary precedence constraints given by a directed acyclic graph. The OPS scheduling problem considered in the present work generalizes problem considered in~\cite{birgin2014milp} in three relevant ways: \textbf{(a)} the concept of precedence among operations of a job is redefined allowing some overlapping; \textbf{(b)} operations are resumable since their processing may be interrupted by periods of unavailability of the machines; and \textbf{(c)} sequence-dependent setup times must be considered. Other than that, operations' release times are considered; and some operations may have pre-defined machines and starting times for execution.

Note that the concept of ``job'' is indirectly addressed in the considered problem. It is assumed that a job is composed of a set of operations and that precedence relations may exist between operations belonging to the same job; while operations belonging to different jobs have no precedence relations among them. Once the set of precedence relations has been defined (as the union of the precedence relations of each individual job), the number~$n$ of jobs plays no explicit role in the problem definition anymore. Let~$o$ and~$m$ be, respectively, the number of operations and machines. Let $V=\{1,2,\dots,o\}$. For each operation~$i \in V$, let $\emptyset \neq F(i) \subseteq \{1,2,...,m\}$ be the subset of machines that can process operation~$i$ and let $p_{ik}$ ($i \in V$, $k \in F(i)$) be the corresponding processing times. Furthermore, let~$A$ be a set of operations pairs $(i,j)$ representing precedence relations. (Note that, by definition, directed arcs in~$A$ correspond to independent directed acyclic graphs (DAG), each one representing the precedence relations among the operations of a job.) Each operation~$i \in V$ is associated with a quantity~$\theta_i \in (0,1]$ that represents the fraction of operation~$i$ that must be completed before starting the execution of any operation~$j$ such that $(i,j) \in A$. In addition, the existence of the pair~$(i,j)$ in~$A$ implies that operation~$j$ can not be completed before the completion of operation~$i$. For each~$i \in V$ and $k \in F(i)$, let $\gamma^F_{ik}$ be the setup time for the case in which operation~$i$ is assigned to machine~$k$ and it is the first operation to be executed (the supra-index ``F'' stands for \textit{first}); and, for each pair of operations $i$ and $j$ and $k \in F(i) \cap F(j)$, let~$\gamma^I_{ijk}$ be the setup time for the case in which both operations are assigned to machine~$k$ with operation~$i$ being the immediate predecessor of operation~$j$ (the supra-index ``I'' stands for \textit{in between}). For each machine~$k$ ($k=1,\dots,m$), let $q_k$ be the number of periods of unavailability, given by $[\underline{u}^k_1,\bar u^k_1], \dots, [\underline{u}^k_{q_k},\bar u^k_{q_k}]$. For each operation~$i \in V$, let~$r_i$ be its release time. Finally, let $T \subseteq V$ be the set of operations with pre-defined starting times $\bar s_i$ for all $i \in T$. The pre-assignment of an operation~$i \in T$ to a machine~$\bar k_i$ does not deserve a special treatment since it can be expressed by simply defining $F(i) = \{ \bar k_i \}$ for all $i \in T$, i.e.,\ the set $F(i)$ of machines that can process an operation~$i \in T$ being a singleton with~$\bar k_i$ as its sole element. The problem consists of assigning each operation~$i \in V$ to a machine $k \in F(i)$ and determining its starting processing time~$s_i$ to satisfy all the problem constraints. A machine can not process more than one operation at a time. The processing of an operation can only be interrupted by periods of unavailability of the machine to which it was assigned. A setup operation can not be interrupted and its completion must coincide with the beginning of the associated operation. The objective is to minimize the makespan. 
  
In addition to the real problem constraints, three modeling  assumptions are being done. In \textbf{Assumption~A1}, it is assumed that the starting time of an operation can not coincide with the beginning nor belong to the interior of a period of unavailability of the machine to which the operation has been assigned. Note that the situations being avoided correspond to the case in which the \textit{actual} operation's starting time coincides with the end of the unavailability period; so all the forbidden cases can be represented by the latter one without loss of generality. In an analogous way, in \textbf{Assumption~A2}, it is assumed that the completion time of an operation can no coincide with the end nor belong to the interior of a period of unavailability of the machine to which the operation has been assigned. The situations being avoided correspond to the case in which the \textit{actual} operation's completion coincides with the beginning of the unavailability period; and the forbidden situations can be represented by the latter one without loss of generality. The third modeling  assumption serves to the purpose of modeling  Assumptions~A1 and~A2. \textbf{Assumption~A3}, in complete accordance with practice, says that all constants that define the problem (processing times, beginning and end of periods of unavailability, starting times of fixed operations, etc.) are rational numbers. Moreover, without loss of generality, by a change of scale or units of measure, it can be assumed that all these constants are integer numbers. With the help of Assumption~A3, the constraint in Assumption~A1 can be modeled requesting that between the starting time of an operation and the beginning of a forthcoming period of unavailability there must be at least a unit of time. In an analogous way, with the help of Assumption~A3, the constraint in Assumption~A2 can be modeled requesting that between the end of a period of unavailability and a forthcoming completion time (of an operation that started to be executed before the period of unavailability) there must be at least a unit of time.

Ignoring the three relevant generalizations listed at the beginning of the current section and named \textbf{(a)}, \textbf{(b)}, and \textbf{(c)}, the OPS scheduling problem can be modeled as a mixed integer linear programming with minor modifications to the MILP model for the FJS scheduling problem with sequence flexibility in~\cite{birgin2014milp}. The model considers binary variables $x_{ik}$ ($i \in V$, $k \in F(i)$) to indicate whether operation~$i$ is assigned to machine~$k$ ($x_{ik} = 1$) or not ($x_{ik} = 0$). It also considers binary variables $y_{ij}$ ($i,j \in V$, $i \neq j$, $F(i) \cap F(j) \neq \emptyset$) to indicate, whenever two operations are assigned to the same machine, which one is processed first. (If~$i$ is processed before~$j$, we have $y_{ij}=1$ and $y_{ji}=0$; while if~$j$ is processed before~$i$, we have $y_{ji}=1$ and $y_{ij}=0$.) Finally, the model uses real variables $s_i$ ($i \in V$) to denote the starting time of operation~$i$ and a variable $C_{\max}$ to represent the makespan. With these variables, the modified MILP model can be written as follows:
\begin{align}
    \mbox{Minimize } C_{\max}    \quad & \quad \label{eq:min_cmax}\\
    \mbox{subject to}            \quad & \quad \nonumber  \\
    \sum_{k \in F(i)} x_{ik} = 1 \quad & \quad i \in V, \label{eq:birgin2}\\ 
    p'_i = \sum_{k \in F(i)} x_{ik} p_{ik} \quad & \quad i \in V, \label{eq:birgin3}\\ 
    s_i + p'_i \leq C_{\max}     \quad & \quad i \in V, \label{eq:birgin1}\\ 
    s_i + p'_i \leq s_j \quad & \quad (i,j) \in A, \label{eq:birgin4}\\  
    y_{ij} + y_{ji} \geq x_{ik} + x_{jk} - 1 \quad & \quad i,j \in V, \; i \neq j, \; k \in F(i) \cap F(j), \label{eq:birgin5}\\[1mm]  
    s_i + p'_i - (1 - y_{ij}) M \leq s_j \quad & \quad i,j \in V, \; i \neq j, \; F(i) \cap F(j) \neq \emptyset, \label{eq:6}\\ 
    s_i \geq r_i                 \quad & \quad i \in V, \label{eq:7}\\
    s_i = \bar s_i               \quad & \quad i \in T. \label{eq:7nova}
\end{align}

The objective function~(\ref{eq:min_cmax}) minimizes the makespan~$C_{\max}$. Constraints~(\ref{eq:birgin2}) say that each operation must be assigned to exactly one machine; and constraints~(\ref{eq:birgin3}) define~$p'_i$ as the processing time of operation~$i$ (that depends on the machine to which it was assigned.) Constraints~(\ref{eq:birgin1}) ensure that every operation is not completed later than~$C_{\max}$. Constraints~(\ref{eq:birgin4}) represent the classical meaning of a precedence constraint, saying that, independently of the machine by which they are being processed, if operation~$i$ precedes operation~$j$, then~$i$ must be completed before~$j$ can be started. (The meaning of ``precedence'' will be redefined in Section~\ref{redefinitionofprecedence}.) Constraints~(\ref{eq:birgin5}) and~(\ref{eq:6}) deal with pairs of operations~$i$ and~$j$ that are assigned to the same machine~$k$. On the one hand, constraints~(\ref{eq:birgin5}) say that if~$i$ and~$j$ were assigned to the same machine~$k$, then~$i$ must be processed before~$j$ or~$j$ must be processed before~$i$, i.e.,\ an order must be established. On the other hand, constraints~(\ref{eq:6}) say that if~$i$ and~$j$ were assigned to the same machine~$k$ and~$i$ is processed before~$j$, then between the starting times of~$i$ and~$j$ there must be enough time to process~$i$. Constraints~(\ref{eq:7}) say that each operation~$i$ can not start to be processed before its release time; and constraints~(\ref{eq:7nova}) fix the start time of the operations whose start times are already pre-determined. 

A summary of the constants and sets that define an instance of the problem as well as of the model's variables can be found in Table~\ref{tab:MILPdescr}. The table can be use as a reference for the forthcoming sections as well.

\begin{table}[ht!]
\caption{Description of the constants, sets, and variables in the MILP formulations of the OPS scheduling problem. Section~\ref{sec:mip} corresponds to the basic model. In Section~\ref{redefinitionofprecedence} the partial overlapping feature is included. Machines' periods of unavailability are included in Section~\ref{sec:unavs} and sequence-dependent setup times are included in Section~\ref{sec:mipsetup}.}
\label{tab:MILPdescr}
\begin{center}
{\scriptsize
\begin{tabular}{|c|m{7cm}|m{7.2cm}|}
\hline
\multicolumn{3}{|c|}{Problem data (constants and sets)}\\
\hline
\hline
\multirow{9}{*}{\S2}&$o$ & number of operations\\
\arrayrulecolor{lightgray}
\cline{2-3}
&$m$ & number of machines\\
\cline{2-3}
&$V=\{1,2,\dots,o\}$ & set of operations\\
\cline{2-3}
&$T \subseteq V$ & set of fixed operations\\
\cline{2-3}
&$F(i) \subseteq \{1,2,\dots,m\}$, $i \in V$ & set of machines that can process operation $i$\\
\cline{2-3}
&$A \subseteq V \times V$ & set of operations' precedence relations \\
\cline{2-3}
&$\bar s_i$, $i \in T$ & starting time for the fixed operation $i$\\
\cline{2-3}
&$r_i$, $i \in V$ & release time of operation $i$\\
\cline{2-3}
&$p_{ik}$, $i \in V$, $k \in F(i)$ & processing time of operation $i$ on machine $k$\\
\hline
\S2.1&$\theta_i \in (0,1]$, $i \in V$ & fraction of operation $i$ that must be processed before an operation $j$ can start to be processed if $(i,j) \in A$\\
\hline
\multirow{2}{*}{\S2.2}&$q_k$, $k \in \{1,2,\dots,m\}$ & number of periods of unavailability of machine $k$\\
\cline{2-3}
&$\underline{u}_{\ell}^k$ and $\bar u_{\ell}^k$, $\ell \in \{1,\dots,q_k\}$ & begin and end of the $\ell$th unavailability period of machine $k$\\
\hline
\multirow{3}{*}{\S2.3}&$\gamma_{ijk}^I$, $i,j \in V$, $i \neq j$, $k \in F(i) \cap F(j)$ & setup time to process operation $j$ right after operation $i$ on machine $k$\\
\cline{2-3}
&$\gamma_{ik}^F$, $i \in V$, $k \in F(i)$ & setup time to process operation $i$ as the first operation to be processed on machine $k$\\
\cline{2-3}
&$B_k = \{ i \in V | k \in F(i) \}$, $k \in \{1,2,\dots,m\}$ & set of operations that can be processed on machine $k$ (auxiliar set that simplifies the presentation of some constraints)\\
\arrayrulecolor{black}
\hline
\multicolumn{3}{c}{}\\
\hline
\multicolumn{3}{|c|}{MILP models' variables}\\
\hline
\hline
\multirow{5}{*}{\S2}&$x_{ik} \in \{0,1\}$, $i \in V$, $k \in F(i)$ & $x_{ik}=1$ if and only if operation $i$ is assigned to machine $k$\\
\arrayrulecolor{lightgray}
\cline{2-3}
&$y_{ij} \in \{0,1\}$, $i,j \in V$, $i \neq j$, $F(i) \cap F(j) \neq \emptyset$ & $y_{ij}=1$ if operations $i$ and $j$ are assigned to the same machine and $i$ is processed before $j$\\
\cline{2-3}
&$p'_i$, $i \in V$ & processing time of operation $i$ (its value depends on the machine to which operation $i$ is assigned)\\
\cline{2-3}
&$s_i$, $i \in V$ & starting time of operation $i$ (it must coincide with $\bar s_i$ if $i \in T$)\\
\cline{2-3}
&$C_{\max}$ & makespan\\
\hline
\S2.1&$\bar p'_i$, $i \in V$ & processing time of operation $i$ that must be completed before an operation $j$ can start to be processed if $(i,j) \in A$ (its value coincides with $\lceil \theta_i p'_i \rceil$)\\
\hline
\multirow{5}{*}{\S2.2}&$c_i$, $i \in V$ & completion time of operation $i$\\
\cline{2-3}
&$\bar c_i$, $i \in V$ & instant at which $\bar p'_i$ units of time of operation $i$ has been processed\\
\cline{2-3}
&$v_{ik\ell}$, $w_{ik\ell}$, $\bar w_{ik\ell} \in \{0,1\}$, $i \in V$, $k \in F(i)$, $\ell \in \{1,\dots,q_k\}$ & auxiliar variables to model constraints on $s_i$, $c_i$, and $\bar c_i$ with respect to the machines' unavailability periods\\
\cline{2-3}
&$u_i$, $i \in V$ & amount of time the machine $k(i)$ to which operation $i$ is assigned is unavailable between $s_i$ and $c_i$\\
\cline{2-3}
&$\bar u_i$, $i \in V$ & amount of time the machine $k(i)$ to which operation $i$ is assigned is unavailable between $s_i$ and $\bar c_i$\\
\hline
\multirow{3}{*}{\S2.3}&$y_{ijk} \in \{0,1\}$, $i,j \in V$, $i \neq j$, $k \in F(i) \cap F(j)$ & $y_{ijk}=1$ if and only if operations $i$ and $j$ are assigned to machine $k$ and $i$ is the immediate predecessor of $j$ (these variables substitute variables $y_{ij}$ from \S2)\\
\cline{2-3}
&$\hat \xi_{ik}$, $\bar \xi_{ik}$, $i \in V$, $k \in F(j)$ & auxiliar variables to model the sequence-dependent setup time\\
\cline{2-3}
&$\xi_i$, $i \in V$ & sequence-dependent setup time of operation $i$\\
\arrayrulecolor{black}
\hline
\end{tabular}}
\end{center}
\end{table}

\subsection{Modeling partial overlap} \label{redefinitionofprecedence}
The simplest feature to be added to model~(\ref{eq:min_cmax}--\ref{eq:7nova}) is the redefinition of precedence. The modification consists in substituting~(\ref{eq:birgin4}) by
\begin{equation} \label{eq:birgin4a}
\bar p'_i = \sum_{k \in F(i)} x_{ik} \lceil \theta_i p_{ik} \rceil, \quad (i,j) \in A,
\end{equation}
\begin{equation} \label{eq:birgin4b}
s_i + \bar p'_i \leq s_j, \quad (i,j) \in A,
\end{equation}
and
\begin{equation} \label{eq:birgin4c}
s_i + p'_i \leq s_j + p'_j, \quad (i,j) \in A.
\end{equation}
Constraints~(\ref{eq:birgin4a}) define new variables $\bar p'_i$ ($i \in V$) whose values coincide, by definition, with $\lceil \theta_i p'_i \rceil$. Constraints~(\ref{eq:birgin4b}) say that, if there exists a precedence constraint saying that operation~$i$ must precede operation~$j$, it means that operation~$j$ can not start before~$100\% \times \lceil \theta_i p'_i \rceil / p'_i$ of operation~$i$ is completed; while constraints~(\ref{eq:birgin4c}) say that operation~$j$ can not be completed before operation~$i$ is completed as well.

\subsection{Modeling machines' unavailabilities} \label{sec:unavs}

Modeling the machines' unavailability starts by considering real variables~$c_i$ and~$\bar c_i$ ($i \in V$) and binary variables~$v_{ik\ell}$, $w_{ik\ell}$, and~$\bar w_{ik\ell}$ ($i \in V$, $k \in F(i)$, $\ell=1,\dots,q_k$). The meaning of variables~$c_i$ and $\bar c_i$ ($i \in V$) follows:
\begin{itemize}
    \item Variable~$c_i$ represents the completion time of each operation~$i$.
    \item Variable~$\bar c_i$ represents the completion time of the fraction of each operation~$i$ that is required to be completed before the execution of an operation~$j$ if there exists~$(i,j) \in A$; i.e.,\ while~$c_i$ represents the completion time of the execution of~$p'_i$ units of time of operation~$i$, $\bar c_i$ represents the completion time of the execution of~$\bar p'_i$ units of time of operation~$i$.
\end{itemize} 
Variables~$v_{ik\ell}$, $w_{ik\ell}$, and~$\bar w_{ik\ell}$ relate~$s_i$, $c_i$, and~$\bar c_i$ with the periods of unavailability of the machines, respectively. The relation we seek between $s_i$ and $v_{ik\ell}$ is ``$v_{ik\ell}=1$ if and only if operation~$i$ is assigned to machine~$k$ and the $\ell$-th period of unavailability of machine~$k$ is to the left of~$s_i$.'' The relation we seek between $c_i$ and $w_{ik\ell}$ is analogous, namely, ``$w_{ik\ell}=1$ if and only if operation~$i$ is assigned to machine~$k$ and the $\ell$-th period of unavailability of machine~$k$ is to the left of~$c_i$''. Finally, the relation we seek between $\bar c_i$ and $\bar w_{ik\ell}$ is ``$\bar w_{ik\ell}=1$ if and only if operation~$i$ is assigned to machine~$k$ and the $\ell$-th period of unavailability of machine~$k$ is to the left of~$\bar c_i$''. These relations can be modeled as follows:
\begin{equation}\label{sc}
    \left.
    \begin{array}{rcl}
    v_{ik\ell} &\leq& x_{ik}\\
    s_i &\leq& \underline{u}^k_{\ell} - 1 + M v_{ik\ell} + M ( 1 - x_{ik} )\\
    s_i &\geq& \bar u^k_{\ell} - M ( 1 - v_{ik\ell} ) - M ( 1 - x_{ik} ) \\[2mm]
    w_{ik\ell} &\leq& x_{ik}\\
    c_i &\leq& \underline{u}^k_{\ell} + M w_{ik\ell} + M ( 1 - x_{ik} )\\
    c_i &\geq& \bar u^k_{\ell} + 1 - M ( 1 - w_{ik\ell} ) - M ( 1 - x_{ik} ) \\[2mm]
    \bar w_{ik\ell} &\leq& x_{ik}\\
    \bar c_i &\leq& \underline{u}^k_{\ell} + M \bar w_{ik\ell} + M ( 1 - x_{ik} )\\
    \bar c_i &\geq& \bar u^k_{\ell} + 1 - M ( 1 - \bar w_{ik\ell} ) - M ( 1 - x_{ik} )
    \end{array}
    \right\}
    i \in V, \; k \in F(i), \; \ell=1,\dots,q_k.
\end{equation}
Note that, in addition, constraints~(\ref{sc}) avoid starting and completion times of an operation to belong to the interior of a period of unavailability of the machine to which the operation was assigned. In fact, constraints~(\ref{sc}) also avoid a starting time to coincide with the beginning of a period of unavailability as well as avoid completion times to coincide with the end of a period of unavailability. Now, we have that
\begin{equation}\label{newui}
u_i = \sum_{k \in F(i)} \sum_{\ell=1}^{q_k} ( w_{ik\ell} - v_{ik\ell} ) (\bar u^k_{\ell} - \underline{u}^k_{\ell}), 
\quad i \in V,
\end{equation}
represents the sum of the unavailabilities between instants~$s_i$ and~$c_i$ (assuming $s_i \leq c_i$) of the machine that process operation~$i$; and
\begin{equation}\label{newui2}
\bar u_i = \sum_{k \in F(i)} \sum_{\ell=1}^{q_k} ( \bar w_{ik\ell} - v_{ik\ell} ) (\bar u^k_{\ell} - \underline{u}^k_{\ell}), 
\quad i \in V,
\end{equation}
represents the sum of the unavailabilities between instants~$s_i$ and~$\bar c_i$ (assuming $s_i \leq \bar c_i$) of the machine that process operation~$i$. 

Now, it is easy to see that the relations
\begin{equation} \label{screl}
\left.
\begin{array}{r}
s_i \leq \bar c_i \leq c_i \\[2mm]
s_i + p'_i + u_i = c_i \\ [2mm]
s_i + \bar p'_i + \bar u_i = \bar c_i
\end{array}
\right\} i \in V,
\end{equation}
must hold. (The two equality constraints in~(\ref{screl}) make clear that the new variables~$c_i$ and~$\bar c_i$ are being included into the model with the only purpose of simplifying the presentation, since every appearance of~$c_i$ and~$\bar c_i$ could be replaced by the left-hand side of the respective equality. In fact, the same remark applies to other equality constraints in the model as, for example, constraints~(\ref{eq:birgin3}) that define $p'_i$, constraints~(\ref{eq:birgin4a}) that define $\bar p'_i$, constraints~(\ref{newui}) that define~$u_i$, and constraints~(\ref{newui2}) that define~$\bar u_i$.) With the inclusion of the new variables~$c_i$ and $\bar c_i$, constraints~(\ref{eq:birgin1}), (\ref{eq:6}), (\ref{eq:birgin4b}), and~(\ref{eq:birgin4c}) can be restated, respectively, as follows:
\begin{equation}\label{new1}
c_i \leq C_{\max}, \quad i \in V, 
\end{equation}
\begin{equation}\label{new3}
c_i - ( 1 - y_{ij} ) M \leq s_j, \quad i,j \in V, \; i \neq j, \; F(i) \cap F(j) \neq \emptyset,
\end{equation}
\begin{equation}\label{new2}
\bar c_i \leq s_j, \quad (i,j) \in A,
\end{equation}
\begin{equation}\label{new2b}
c_i \leq c_j, \quad (i,j) \in A.
\end{equation}
    
Summing up, with new real variables~$\bar p'_i$, $c_i$, and $\bar c_i$ ($i \in V$) and new binary variables~$v_{ik\ell}$, $w_{ik\ell}$, and~$\bar w_{ik\ell}$ ($i \in V$, $k \in F(i)$, $\ell=1,\dots,q_k$), the model that includes the machines' unavailable periods consists in minimizing~(\ref{eq:min_cmax}) subject to (\ref{eq:birgin2},\ref{eq:birgin3},\ref{eq:birgin5},\ref{eq:7},\ref{eq:7nova},\ref{eq:birgin4a},\ref{sc},\ref{newui},\ref{newui2},\ref{screl},\ref{new1},\ref{new3},\ref{new2},\ref{new2b}).

\subsection{Modeling sequence-dependent setup time} \label{sec:mipsetup}

In order to model the setup, binary variables $y_{ij}$ will be replaced by binary variables $y_{ijk}$ ($i,j \in V$, $i \neq j$, $k \in F(i) \cap F(j)$). The idea is that $y_{ijk}=1$ if and only if operations~$i$ and~$j$ are assigned to machine~$k$ and~$i$ is the immediate predecessor of~$j$. If, for each machine~$k$, we define the set $B_k = \{ i \in V \; | \; k \in F(i) \}$, the required characterization of variables~$y_{ijk}$ can be achieved with the following constraints:
\begin{equation} \label{setup1}
\left.
\begin{array}{rcl}
y_{ijk} &\leq& x_{ik} \\
y_{ijk} &\leq& x_{jk} \\
\end{array}
\right\}
\quad k=1,\dots,m, \; i,j \in B_k, \; i \neq j,
\end{equation}
\begin{equation} \label{setup2}
\sum_{i,j \in B_k, i \neq j} y_{ijk} = \sum_{i \in B_k} x_{ik} - 1, \quad k=1,\dots,m,
\end{equation}
\begin{equation} \label{setup3}
\sum_{j \in B_k, j \neq i} y_{ijk} \leq 1, \quad k=1,\dots,m, \; i \in B_k,
\end{equation}
\begin{equation} \label{setup4}
\sum_{i \in B_k, i \neq j} y_{ijk} \leq 1, \quad k=1,\dots,m, \; j \in B_k.
\end{equation}
Constraints~(\ref{setup1}) say that, if operation~$i$ or operation~$j$ is not assigned to machine~$k$, then $y_{ijk}$ must be equal to zero. For any machine~$k$, $\sum_{i \in B_k} x_{ik}$ represents the number of operations assigned to it. Then, constraints~(\ref{setup2}) say that, for every machine, the number of precedence constraints related to the operations assigned to it must be one less than the number of operations assigned to it. Finally, constraints~(\ref{setup3}) and~(\ref{setup4}) say that each operation precedes at most one operation and it is preceded by at most one operation, respectively. In fact, it is expected each operation to precede and to be preceded by exactly one operation unless the first and the last operations in the scheduling of each machine. 

The substitution of variables~$y_{ij}$ ($i,j \in V$, $i\neq j$, $F(i) \cap F(j) \neq \emptyset$) by variables $y_{ijk}$ ($i,j \in V$, $i \neq j$, $k \in F(i) \cap F(j)$) consists in replacing constraints~(\ref{eq:birgin5}) with constraints~(\ref{setup1},\ref{setup2},\ref{setup3},\ref{setup4}) and replacing constraints~(\ref{new3}) with
\begin{equation}\label{new4}
c_i - \left( 1 - \sum_{k \in F(i) \cap F(j)} y_{ijk} \right) M \leq s_j, 
\quad i,j \in V, \; i \neq j, \; F(i) \cap F(j) \neq \emptyset.
\end{equation}
Up to this point, we have a model, alternative to the one we already had, but more suitable to include the setup feature.
    
Now, observe that, if operation~$j$ is assigned to machine~$k$ and 
\[
\sum_{i \in B_k, i \neq j} y_{ijk} = 0,
\]
this means that~$j$ is the first operation to be processed by machine~$k$. In this case, the sequence-dependent setup time required to process operation~$j$ on machine~$k$ is given by~$\gamma^F_{jk}$. On the other hand, if operation~$j$ is assigned to machine~$k$ and 
\[
\sum_{i \in B_k, i \neq j} y_{ijk} = 1,
\]
then operation~$j$ is preceded by the only operation~$i$ such that~$y_{ijk}=1$. In this case, the sequence-dependent setup time required to process operation~$j$ on machine~$k$ is given by
\[
\sum_{i \in B_k, i \neq j} y_{ijk} \gamma^I_{ijk}.
\]
Thus, if we define
\begin{equation} \label{setup5}
\hat \xi_{jk} = \left( \sum_{i \in B_k, i \neq j} y_{ijk} \gamma^I_{ijk} \right) +
\left( 1 - \sum_{i \in B_k, i \neq j} y_{ijk} \right) \gamma^F_{jk}, \quad j \in V, \; k \in F(j),
\end{equation}
we have that, if operation~$j$ is assigned to machine~$k$, then $\hat \xi_{jk}$ corresponds to the sequence-dependent setup time required to process operation~$j$ on machine~$k$. Note that, if operation~$j$ is \textit{not} assigned to machine~$k$, by~(\ref{setup5}), we have $\hat \xi_{jk} = \gamma^F_{jk}$. So, defining new real variables~$\bar \xi_{jk}$ ($j \in V$, $k \in F(j)$) and constraints
\begin{equation} \label{setup6}
\left.
\begin{array}{rcccl}
0 &\leq& \bar \xi_{jk} &\leq& M x_{jk}\\
\hat \xi_{jk} - M ( 1 - x_{jk} ) &\leq& \bar \xi_{jk} &\leq& \hat \xi_{jk}\\
\end{array}
\right\}
j \in V, \; k \in F(j),
\end{equation}
we have that, if operation~$j$ is assigned to machine~$k$, then its sequence-dependent setup time is given by~$\bar \xi_{jk}$; while $\bar \xi_{jk}=0$ if operation~$j$ is not assigned to machine~$k$. Thus,
\begin{equation} \label{setup7}
\xi_j = \sum_{k \in F(j)} \bar \xi_{jk}, \quad j \in V,
\end{equation}
represents the sequence-dependent setup time of operation~$j$.

For including the setup time into the model, we need to replace constraints~(\ref{new4}) with
\begin{equation}\label{newnew4}
c_i - \left( 1 - \sum_{k \in F(i) \cap F(j)} y_{ijk} \right) M \leq s_j - \xi_j, \quad i,j \in V, \; i \neq j, \; F(i) \cap F(j) \neq \emptyset,
\end{equation}
to say that between the end of an operation and the beginning of the next operation there must be enough time to process the corresponding setup. It remains to say that a setup operation can not be interrupted and that its completion time must coincide with the starting time of the operation itself. This requirement corresponds to  constraints
\begin{equation} \label{precisa}    
\left.
\begin{array}{rcl}
s_i - \xi_i &\leq& \underline{u}^k_{\ell} - 1 + M v_{ik\ell} + M ( 1 - x_{ik} )\\
s_i - \xi_i &\geq& \bar u^k_{\ell} - M ( 1 - v_{ik\ell} ) - M ( 1 - x_{ik} ) \\[2mm]
\end{array}
\right\}
i \in V, \; k \in F(i), \; \ell=1,\dots,q_k
\end{equation}
plus
\begin{equation} \label{aquefaltava}
s_i \geq \xi_i, \quad i \in V.
\end{equation}
Note that constraints~(\ref{precisa}), roughly speaking, say that, if~$s_i$ is in between two unavailabilities, then~$s_i - \xi_i$ must be in between the same two unavailabilities. Combining~(\ref{precisa}) with~(\ref{sc}) results in
\begin{equation}\label{scnew}
\left.
\begin{array}{rcl}
v_{ik\ell} &\leq& x_{ik}\\
s_i &\leq& \underline{u}^k_{\ell} - 1 + M v_{ik\ell} + M ( 1 - x_{ik} )\\
s_i  - \xi_i &\geq& \bar u^k_{\ell} - M ( 1 - v_{ik\ell} ) - M ( 1 - x_{ik} ) \\[2mm]
w_{ik\ell} &\leq& x_{ik}\\
c_i &\leq& \underline{u}^k_{\ell} + M w_{ik\ell} + M ( 1 - x_{ik} )\\
c_i &\geq& \bar u^k_{\ell} + 1 - M ( 1 - w_{ik\ell} ) - M ( 1 - x_{ik} ) \\[2mm]
\bar w_{ik\ell} &\leq& x_{ik}\\
\bar c_i &\leq& \underline{u}^k_{\ell} + M \bar w_{ik\ell} + M ( 1 - x_{ik} )\\
\bar c_i &\geq& \bar u^k_{\ell} + 1 - M ( 1 - \bar w_{ik\ell} ) - M ( 1 - x_{ik} )
\end{array}
\right\}
i \in V, \; k \in F(i), \; \ell=1,\dots,q_k.
\end{equation}

Summing up, with new real variables~$\hat \xi_{jk}$, $\bar \xi_{jk}$ ($j \in V$, $k \in F(j)$), and $\xi_j$ ($j \in V$) and new binary variables $y_{ijk}$ ($i,j \in V$, $i\neq j$, $F(i) \cap F(j) \neq \emptyset$) that substitute the binary variables $y_{ij}$ ($i,j \in V$), the updated model that includes the sequence-dependent setup time consists in minimizing~(\ref{eq:min_cmax}) subject to (\ref{eq:birgin2},\ref{eq:birgin3},\ref{eq:7},\ref{eq:7nova},\ref{eq:birgin4a},\ref{newui},\ref{newui2},\ref{screl},\ref{new1},\ref{new2},\ref{new2b},\ref{setup1},\ref{setup2},\ref{setup3},\ref{setup4},\ref{setup5},\ref{setup6},\ref{setup7},\ref{newnew4},\ref{aquefaltava},\ref{scnew}). Con\-straints (\ref{setup6}), (\ref{newnew4}), and~(\ref{scnew}) depend on a ``sufficiently large'' constant~$M$ whose value needs to be defined. In~(\ref{setup6}), a sufficiently large value for $M$ is given by
\[
M_1 = 
\max \left\{ 
\max_{j \in V, \; k \in F(j)} \left\{ \gamma^F_{jk} \right\}, 
\max_{i,j \in V, \; i \neq j, \; F(i) \cap F(j) \neq \emptyset} \left\{ \gamma^I_{ijk} \right\} 
\right\}.
\]
In~(\ref{newnew4}), a sufficiently large value for $M$ is given by any upper bound for the optimal~$C_{\max}$ like, for example, 
\[
M_2 = \max_{\{k=1,\dots,m\}}\{ \bar u^k_{q_k} \} +
\sum_{j \in V} \max_{k \in F(j)} \left\{ p_{jk} + \max \left\{ \gamma^F_{jk}, \max_{\{i \in V | k \in F(i)\}} \left\{ \gamma^I_{ijk} \right\} \right\} \right\}.
\]
(Note that $M_2$ is a loose upper bound for the optimal~$C_{\max}$.) The same value~$M_2$ can be used in the second, the fifth, and the eighth inequalities in~(\ref{scnew}); while, in the third, the sixth, and the ninth inequalities in~(\ref{scnew}), it can be used
\[
M_3 = \max_{\{k=1,\dots,m\}}\{ \bar u^k_{q_k} \}.
\]

\section{Constraint programming formulation}\label{sec:cp}

Constraint Programming~\citep{HandbookCP2006} is a powerful paradigm for solving combinatorial problems; and it is particularly attractive for problems that do not have a simple formulation in terms of linear constraints, as it is the case of the OPS scheduling problem being considered in the present work. CP Optimizer~\citep{Laborie2018} is an optimization engine based on CP that extends classical CP with a few mathematical concepts that make it easier to model scheduling problems while providing an interesting problem's structure to its automatic search algorithm. The automatic search is an exact algorithm (so it produces optimality proofs) that internally uses some metaheuristics, mainly the Self-Adapting Large-Neighborhood Search~\citep{Laborie2007}, to quickly produce good quality solutions that help to prune the search space. A formulation of the OPS scheduling problem suitable to be solved with CP Optimizer is presented below. The CP Optimizer concepts will be briefly described as soon as they appear in the formulation. For more details, please refer to~\cite{Laborie2018} and to the CP Optimizer reference manual.

A formulation using the concepts of CP Optimizer equivalent to the MILP model~(\ref{eq:min_cmax}--\ref{eq:7nova}) follows:
\begin{align}
\mbox{Minimize } \max \limits_{i \in V} \cpoEndOf(o_{i}) \quad & \quad \label{eqcpo:3}\\
\mbox{subject to} & \quad \nonumber\\
\cpoEndBeforeStart(o_{i},o_{j}), & \quad (i,j) \in A,\label{eqcpo:4}\\ 
\cpoAlternative(o_{i},[a_{ik}]_{k \in F(i)}), & \quad i \in V, \label{eqcpo:5}\\ 
\cpoNoOverlap([a_{ik}]_{i \in V : k \in F(i)}), & \quad k=1,\dots,m, \label{eqcpo:6}\\ 
\cpoStartOf(o_{i}) \geq r_i, & \quad i \in V, \label{eqcpo:6b}\\ 
\cpoStartOf(o_{i}) = \bar s_i, & \quad i \in T, \label{eqcpo:6c}\\ 
\cpoInterval o_{i}, & \quad i \in V, \label{eqcpo:1}\\ 
\cpoInterval \ a_{ik}, \ \cpoOptional, \ \cpoSize=p_{ik}, & \quad i \in V, \; k \in F(i). \label{eqcpo:2}  
\end{align}
Decision variables of the problem are described in~(\ref{eqcpo:1}) and~(\ref{eqcpo:2}). In~(\ref{eqcpo:1}), an interval variable $o_i$ for each operation $i$ is defined. In~(\ref{eqcpo:2}), an ``optional'' interval variable $a_{ik}$ is being defined for each possible assignment of operation~$i$ to a machine $k \in F(i)$. Optional means that the interval variable may exist or not; and the remaining of the constraint says that, in case it exists, its size must be $p_{ik}$. The objective function~(\ref{eqcpo:3}) is to minimize the makespan, given by the maximum end value of all the operations represented by the interval variables~$o_i$. Precedence constraints between operations are posted as $\cpoEndBeforeStart$ constraints between interval variables in constraints~(\ref{eqcpo:4}). Constraints~(\ref{eqcpo:5}) state that each operation~$i$ must be allocated to exactly one machine~$k \in F(i)$ that is, one and only one interval variable~$a_{ik}$ must be present and the selected interval~$a_{ik}$ will start and end at the same values as interval~$o_i$. Constraints~(\ref{eqcpo:6}) state that, for a machine $k$, the intervals $a_{ik}$ representing the assignment of the operations to this machine do not overlap. (It should be noted that this $\cpoNoOverlap$ constraints actually create a hidden sequence variable on the intervals $a_{ik}$. More details on sequence variables will be given on Section~\ref{sec:setup}.) Finally, constraints~(\ref{eqcpo:6b}) say that each operation~$i$ can not start to be processed before its release time; and constraints~(\ref{eqcpo:6c}) fix the starting time of the operations whose starting times are already pre-determined. 

The CP Optimizer model~(\ref{eqcpo:3}--\ref{eqcpo:2}) can be strengthened by a redundant constraint stating that, at any moment in time, there are never more than~$m$ machines being used simultaneously. This constraint (that implicitly apply to discrete instants in time) is given by
\begin{equation} \label{eqcpo:7}
\sum_{i \in V} \cpoPulse(o_{i},1) \leq m.
\end{equation}
It is worth noting that, in CP as well as in MILP formulations, adding redundant constraints is a common technique to get a stronger reduction of the domains of variables at each search node, which results in a better pruning of the search space.

\subsection{Modeling machines' unavailabilities}

As the operations are suspended by the unavailabilities of the machines, the definition of the interval variable~$a_{ik}$ must be modified by considering as intensity function a step function~$U_k$ that represents the unavailability of machine~$k$. In CP Optimizer, step functions are constant structures of the model that are represented by a set of steps associated with a value. The value of the step function~$U_k$ is 0\% when machine~$k$ is unavailable, i.e.,\ on time windows $[\underline{u}^k_{1},\bar u^k_{1}],...,[\underline{u}^k_{q_k},\bar u^k_{q_k}]$, and 100\% in between these time windows. So, including the machines' unavailabilities in model~(\ref{eqcpo:3}--\ref{eqcpo:2}) simply consists in replacing~(\ref{eqcpo:2}) with
\begin{equation} \label{eqcpoops:B2}
\cpoInterval \ a_{ik}, \ \cpoOptional, \ \cpoSize=p_{ik}, \ \cpoIntensity = U_k, 
\; i \in V, k \in F(i).
\end{equation}

\subsection{Modeling partial overlap}

\begin{figure}[ht!]
\centering
\footnotesize
\includegraphics[width=\textwidth]{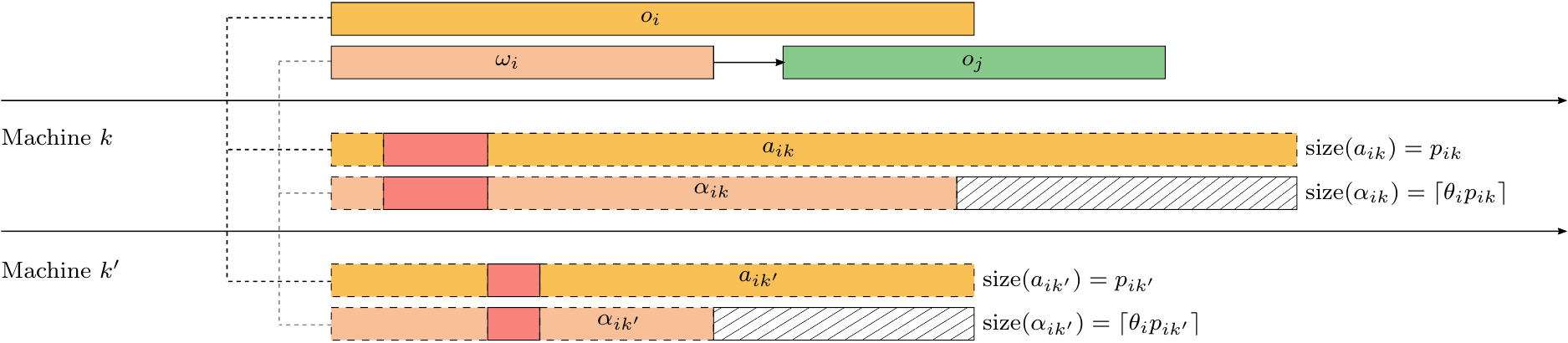} 
\caption{Illustration of the relation between the interval variables~$o_i$, $a_{ik}$, and $\alpha_{ik}$. In the picture, there are two alternative machines~$k$ and~$k'$ to process operation~$i$; and operation~$i$ is assigned to machine~$k'$ since the size of~$\omega_i$ coincides with the size of~$\alpha_{ik'}$. The size of~$\omega_i$ corresponds to the proportion of operation~$i$ that must be processed before operation~$j$ can start to be processed. (The picture assumes that there is a precedence constraint between operations~$i$ and~$j$.)}
\label{fig:cpmodel1}
\end{figure}

An additional feature of the OPS scheduling problem is that operations subject to precedence constraints may partially overlap. Each interval variable~$o_i$ (resp.\ $a_{ik}$) is associated with an additional interval variable~$\omega_i$ (resp.\ $\alpha_{ik}$) that represents the proportion of  operation~$i$ that has to be processed before any of its successors can start. The size of the optional interval variable~$\alpha_{ik}$ is defined as $\lceil \theta_i p_{ik} \rceil$ (see~(\ref{eqcpoops:1}) and~(\ref{eqcpoops:2})); and interval variables~$a_{ik}$ and~$\alpha_{ik}$ have the same presence status (see~(\ref{eqcpoops:3})). Interval~$\omega_i$ is the alternative between all the interval variables~$\alpha_{ik}$ (see~(\ref{eqcpoops:4})); and interval variables~$o_i$ and~$\omega_i$ start at the same time (see~(\ref{eqcpoops:5})). This way, when operation~$i$ is allocated to machine~$k$, both interval variables~$a_{ik}$ and~$\alpha_{ik}$ are present, the size of $\alpha_{ik}$, that represents the proportion of operation~$i$ to be executed before any of its successors, is equal to $\lceil \theta_i p_{ik} \rceil$, and interval variable~$\omega_{i}$ is synchronized with the start and end of~$\alpha_{ik}$. (See Figure~\ref{fig:cpmodel1}.) Precedence constraints are posted between $\omega_i$ and its successors in~(\ref{eqcpoops:6}). Constrains~(\ref{eqcpoops:6bis}) say that operation~$j$ can not be completed before the completion of operation~$i$ if $(i,j) \in A$.
\begin{align}
\cpoInterval \ \omega_{i}, & \quad i \in V, \label{eqcpoops:1}\\ 
\cpoInterval \ \alpha_{ik}, \ \cpoOptional, \ \cpoSize = \lceil \theta_i p_{ik} \rceil, 
\cpoIntensity = U_k, & \quad i \in V, k \in F(i), \label{eqcpoops:2}\\
\cpoPresenceOf(a_{ik}) == \cpoPresenceOf(\alpha_{ik}), & 
\quad i \in V, k \in F(i), \label{eqcpoops:3}\\
\cpoAlternative(\omega_{i},[\alpha_{ik}]_{k \in F(i)}), & \quad i \in V, \label{eqcpoops:4}\\
\cpoStartAtStart(\omega_{i},o_{i}), & \quad i \in V, \label{eqcpoops:5}\\
\cpoEndBeforeStart(\omega_{i},o_{j}), & \quad (i,j) \in A, \label{eqcpoops:6}\\
\cpoEndBeforeEnd(o_{i},o_{j}), & \quad (i,j) \in A. \label{eqcpoops:6bis}
\end{align}

The model that includes the machines unavailabilities as well as the partial overlapping consists in minimizing~(\ref{eqcpo:3}) subject to (\ref{eqcpo:5}--\ref{eqcpo:1}), (\ref{eqcpoops:B2}), and (\ref{eqcpoops:1}--\ref{eqcpoops:6bis}) (that substitute constraints~(\ref{eqcpo:4})); constraints~(\ref{eqcpo:7}) being, as already mentioned, optional.

\subsection{Modeling the sequence-dependent setup time} \label{sec:setup}

\begin{figure}[ht!]
    \centering
    \footnotesize 
    \includegraphics[width=\textwidth]{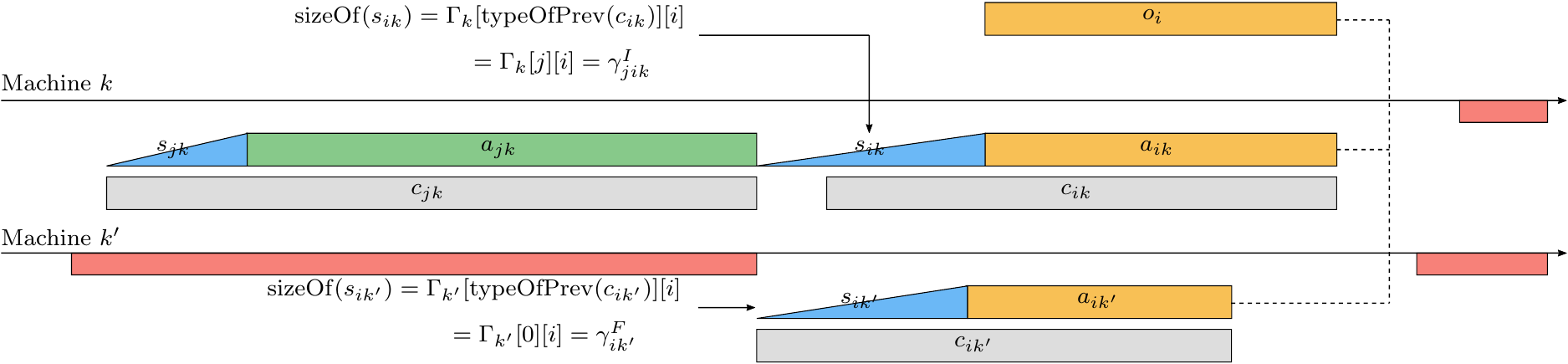} 
    \caption{Illustration of the relation between interval variables~$o_i$, $a_{ik}$, and~$c_{ik}$ in the CP Optimizer model. The illustration also shows the explicit inclusion of the setup interval variables~$s_{ik}$. The right hand side of the graphic shows two alternatives for processing operation~$i$ on machine~$k$ or on machine~$k'$. The left hand side illustrates the computation of the sequence-dependent setup time.}
    \label{fig:cpmodel2}
\end{figure}

Another additional feature of the OPS scheduling problem is the notion of setup time and setup activities between consecutive operations executed on a machine. The usual formulation for sequence-dependent setup times in CP Optimizer is to use a sequence variable that permits to associate an integer type with each interval variable in the sequence and to post a no-overlapping constraint on this sequence with a transition distance matrix. This is the purpose of the sequence variable $SA_k$ defined in~(\ref{eqcpoops:12b}). This sequence variable is defined on all the interval variables~$a_{ik}$ on machine~$k$ (optional interval variables representing the possible execution of operation $i$ on machine~$k$). Each interval $a_{ik}$ is associated with a type~$i$ in the sequence variable. A non-overlapping constraint is posted in~(\ref{eqcpoops:18}), specifying the transition distance matrix $\Gamma^I_k$ defined as $\Gamma^I_k[i][j]=\gamma^I_{ijk}$. These constraints ensure that operations allocated to machine~$k$ do not overlap and that a minimal setup time of $\Gamma^I_k[i][j]$ must elapse between any two consecutive operations~$i$ and~$j$. The definition of the sequence variable and the non-overlapping constraints are given by
\begin{align}
\cpoSequence \ SA_k \ \cpoSeqOn \ [a_{ik}]_{i \in V : k \in F(i)}, \ 
\cpoSeqTypes \ [i]_{i \in V : k \in F(i)}, & \quad k=1,\dots,m, \label{eqcpoops:12b}\\
\cpoNoOverlap(SA_k,\Gamma^I_k), & \quad k=1,\dots,m. \label{eqcpoops:18}
\end{align}
An additional feature of the OPS scheduling problem is that the setup activities are also subject to the unavailability of the machines and, in particular, the setup cannot be interrupted by an unavailability time window. Because of these complex constraints on the setups, the setup operations need to be explicitly represented as interval variables in the model. The explicit representation of the setup operations as interval variables also allows us to model the setup of the first operation being processed by a machine, a case that is not covered by~(\ref{eqcpoops:12b},\ref{eqcpoops:18}). For each interval~$a_{ik}$, we define in~(\ref{eqcpoops:10}) an interval variable $s_{ik}$ that represents the setup activity just before operation $a_{ik}$; and we establish in~(\ref{eqcpoops:17}) that the end of~$s_{ik}$ must coincide with the beginning of~$a_{ik}$. For each interval~$a_{ik}$, we also define in~(\ref{eqcpoops:11}) an interval variable~$c_{ik}$ that covers both~$s_{ik}$ and~$a_{ik}$ (see~(\ref{eqcpoops:13})). Interval variables $a_{ik}$, $s_{ik}$, and $c_{ik}$ have the same presence status (see~(\ref{eqcpoops:14}) and~(\ref{eqcpoops:15})). A sequence variable $SC_k$ representing the non-overlapping sequence of intervals $c_{ik}$ on machine $k$ is defined in~(\ref{eqcpoops:12}), the non-overlapping being represented in~(\ref{eqcpoops:19}). The size of the setup activity~$s_{ik}$ depends both on the type of the previous operation on sequence~$SC_k$ and the type of operation~$i$; and its value is given in~(\ref{eqcpoops:16}) by a matrix~$\Gamma_k$. Matrix~$\Gamma_k$ is a matrix with row index starting from~$0$ and defined as the concatenation of matrices~$\Gamma^F_k$ and~$\Gamma^I_k$, i.e.,\ $\Gamma_k$ consists in~$\Gamma^I_k$ with an additional $0$-th row given by $\Gamma^F_k[i] = \gamma^F_{ik}$. By convention, when operation~$i$ is assigned to a machine~$k$ and it is the first activity executed by the machine, the type of the previous operation on sequence~$SC_k$ is $0$ so that the size of the setup activity $s_{ik}$ is $\Gamma_k[0][i] = \Gamma^F_k[i] = \gamma^F_{ik}$. Figure~\ref{fig:cpmodel2} illustrates the extent of the setup time concerning the predecessor type.

Finally, a set of constraints is added to ensure the behavior of operations and setup activities with respect to machine unavailable periods. The intensity function $U_k$ used in the definition of the interval variables representing operations (\ref{eqcpoops:B2}) states that operations are suspended by unavailabilities. Additionally, constraints~(\ref{eqcpoops:21}) and~(\ref{eqcpoops:22}) establish that operations cannot start or end during an unavailability period; whereas constraints~(\ref{eqcpoops:20}) say that setup activities cannot overlap unavailability periods. The new set of constraints follows:
\begin{align}
\cpoInterval \ s_{ik}, \ \cpoOptional, & \quad i \in V, k \in F(i), \label{eqcpoops:10}\\
\cpoInterval \ c_{ik}, \ \cpoOptional, & \quad i \in V, k \in F(i), \label{eqcpoops:11}\\
\cpoSequence \ SC_k \ \cpoSeqOn \ [c_{ik}]_{i \in V: k \in F(i)}, \ 
\cpoSeqTypes \ [i]_{i \in V: k \in F(i)}, & \quad k=1,\dots,m, \label{eqcpoops:12} \\
\cpoSpan(c_{ik},[s_{ik},a_{ik}]), & \quad i \in V, k \in F(i), \label{eqcpoops:13}\\
\cpoPresenceOf(a_{ik})==\cpoPresenceOf(c_{ik}), & \quad i \in V, k \in F(i), \label{eqcpoops:14}\\
\cpoPresenceOf(a_{ik})==\cpoPresenceOf(s_{ik}), & \quad i \in V, k \in F(i), \label{eqcpoops:15}\\
\cpoSizeOf(s_{ik})==\Gamma_k[\cpoTypeOfPrev(SC_k,c_{ik})][i], & \quad i \in V, k \in F(i), \label{eqcpoops:16}\\
\cpoEndAtStart(s_{ik},a_{ik}), & \quad i \in V, k \in F(i), \label{eqcpoops:17}\\
\cpoNoOverlap(SC_k), & \quad k=1,\dots,m, \label{eqcpoops:19}\\
\cpoForbidStart(a_{ik},U_k), & \quad i \in V, k \in F(i), \label{eqcpoops:21}\\
\cpoForbidEnd(a_{ik},U_k), & \quad i \in V, k \in F(i), \label{eqcpoops:22}\\
\cpoForbidExtent(s_{ik},U_k), & \quad i \in V, k \in F(i). \label{eqcpoops:20}
\end{align}

Summing up, the full CP Optimizer formulation of the OPS scheduling problem is given by the minimization of~(\ref{eqcpo:3}) subject to constraints~(\ref{eqcpo:5}--\ref{eqcpo:1},\ref{eqcpoops:B2},\ref{eqcpoops:1}--\ref{eqcpoops:6bis},\ref{eqcpoops:10}--\ref{eqcpoops:20}); constraints~(\ref{eqcpo:7}) being, as already mentioned, optional. In the same sense, constraints~(\ref{eqcpoops:12b},\ref{eqcpoops:18}), superseded by constraints~(\ref{eqcpoops:10}--\ref{eqcpoops:20}), can be considered optional. It is expected that keeping them would result in a stronger inference in the resolution process due to the direct formulation of the minimal distance~$\Gamma^I_k$ between consecutive operations on each machine~$k$.

\section{Discussion}\label{sec:disc}

In this section, we elaborate on the relation between the main components of the MILP and the CP Optimizer formulations of the OPS scheduling problem. The main difference between both formulations is that, in the CP Optimizer model, the number of explicit variables and constraints is $O(o \, m)$; while, in the MILP formulation, the number of variables and constraints is, in the worst case, $O(o^2 m + o \sum_{k=1}^m q_k)$, where $o$ is the number of operations, $m$ is the number of machines, and $q_k$ is the number of periods of unavailability of machine~$k$. A tighter bound is given by $O(|A| + \sum_{k=1}^m ( |B_k|^2 + |B_k| q_k ) )$, where $A$ is the set of precedence relations and, for each machine~$k$, $B_k$ is the set of operations that could be processed by it. On the other hand, every constraint in the MILP model involves a constant number of variables or a number that is, in the worst case, $O(o + \sum_{k=1}^m q_k)$; while, in the CP Optimizer model, each non-overlap constraint on the sequence $SA_k$ ($k=1,\dots,m$) involves a dense matrix $\Gamma^I_k$ of size~$o^2$. So, as expected, the non-overlapping constraints in the CP Optimizer model involve, as a whole, dealing with the $O(o^2 m)$ given setup times. (The same is true, in a similar way, for the sequence variables $SC_k$.) Another difference between the models is that the CP Optimizer model strongly relies on the integrality of all the constants that define an instance, i.e.,\ it assumes that the processing times of the operation on the machines, the fixed starting times, the setup times, the beginning and the end of the machines' unavailability periods, and the release times are all integer values. If, on the one hand, this is \textit{not} a requirement of the MILP formulation; on the other hand, this is a usual assumption that can be accomplished, in practice, by changing the constants' unit of measure. The same is not true in the required interpretation of partial overlapping. In the considered definition of overlapping, the interpretation of constant~$\theta_i$ is that if~$p'_i$ is the processing time of operation~$i$ on the machine to which it was assigned, then operation~$i$ must be processed at least $\lceil \theta_i p'_i\rceil$ units of time before any successor can start to be processed. The fact of using $\lceil \cdot \rceil$ in the interpretation of overlapping is needed in the CP Optimizer formulation, due to the integrality assumption; while it is not relevant at all in the MILP formulation.

One of the most relevant components of both models is the one that represents the duration of each operation, where by duration we mean its processing time on the machine to which it was assigned plus the duration of the unavailabilities of such machine that interrupt its execution. In the CP Optimizer model, this object is represented by the interval variables~$o_i$. The equivalent object in the MILP formulation is given by the starting time~$s_i$ and the completion time~$c_i$, that correspond, respectively, to the beginning and the end of the interval variable $o_i$ of the CP Optimizer model. In the MILP model, the time elapsed between~$s_i$ and~$c_i$ is divided into the processing time of the operation itself, represented by~$p'_i$, and the sum of the unavailability windows in between~$s_i$ and~$c_i$, represented by $u_i$. In the CP Optimizer model, this distinction is made with the help of the indicator function~$U_k$. 

Another object present in both formulations and strongly related to the duration of an operation is the duration of the fraction of an operation that must be processed before any successor of it can start to be processed. In the CP Optimizer model, this role is played by the interval variables~$\omega_i$. In the MILP model, the beginning and the end of this interval are given by variables~$s_i$ and~$\bar c_i$, respectively. Not by coincidence, there is a constraint in the CP Optimizer model saying that the beginning of $o_i$ and~$\omega_i$ must coincide. (In the MILP model the beginning of these two intervals is represented by the same variable~$s_i$.) Once again, the time elapsed between~$s_i$ and~$\bar c_i$ is divided in the MILP model in the actual processing of the operation (represented by~$\bar p'_i$) and the aggregated time of the unavailability windows that interrupt the execution of the operation, represented by~$\bar u_i$). The same effect is, once again, obtained in the CP Optimizer model through the usage of the indicator function~$U_k$. 

The assignment of an operation to a machine is modeled in the MILP formulation with binary variables~$x_{ik}$. Then, with the same variables, the machine-dependent processing time $p_i$ of operation~$i$ is established. The same thing is done in the CP Optimizer formulation with the ``alternative'' constraint that relates the interval variables $o_i$ with one and only one variable $a_{ik}$ for some $k \in F(i)$, whose duration is~$p_{ik}$. Note that, in the MILP model, the $x_{ik}$ variables also only exist if $k \in F(i)$. Both situations correspond to the disjunction that is modeled in a classical way in the MILP model; while it corresponds to a primitive of the CP Optimizer modeling language in the other case. The non-overlapping between operations being processed by the same machine~$k$ is once again modeled with the help of binary variables $y_{ijk}$ in the MILP model. In CP Optimizer, the ``noOverlap'' constraints are posted on a hidden ``sequence'' variable that is defined over a set of interval variables, each interval variable in the set being associated with an integer ``type.'' The ``noOverlap'' constraints state that the intervals of the sequence must be ordered as a set of non-overlapping intervals (that typically represent the different operations on a machine). The ``type'' of the types of the intervals may, for instance, be used to index a setup time matrix that represents a minimal distance between consecutive intervals in the sequence. In the MILP model, the setup feature is achieved with the help of the $y_{ijk}$ variables. All other elements of both models also relate similarly.

As it was described, the MILP and the CP Optimizer model represent the OPS scheduling problem using equivalent structures; the difference lies on the fact that the CP Optimizer modeling language provides to the CP Optimizer solver a more compact model and much more structure of the problem than the one that the MILP model could give to a general-purpose MILP solver.

\section{Numerical experiments}\label{sec:exp}

The numerical experiments in the present section have three goals. Experiments with small-sized instances aim to compare the efficiency and effectiveness of the commercial solvers IBM ILOG CPLEX and IBM ILOG CP Optimizer (both included in version 12.9), when applied to instances modeled with the Python API DOcplex 2.10.155 library. Experiments with medium-sized instances aim to determine the size of the instances for which optimality can be proved with the commercial solvers IBM ILOG CP Optimizer. Experiments with large-sized instances, that are of the size of real instances of the OPS scheduling problem, are also considered. The goal of these experiments is to evaluate the possibility of using feasible solutions found by the IBM ILOG CP Optimizer in practice. Section~\ref{secinst} describes the generation of instances; while Section~\ref{secres} presents the numerical results with small-, medium, and large-sized instances. The implementation of the MILP and the CP Optimizer formulations, the generator of random instances, and all the generated instances are freely available at \url{https://willtl.github.io/ops}.

\subsection{Generation of instances}
\label{secinst}

Numerical experiments with the introduced models were performed on random instances. All constants that define an instance are numbers randomly chosen with uniform distribution in a predefined interval. So, from now on, whenever ``chosen'', ``random'' or ``randomly chosen'' is written, it should be read ``randomly chosen with uniform distribution''. It should be noted that, although random, instances are generated in such a way they preserve the characteristics of the real instances of the OPS scheduling problem. Three different sets with small-, medium-, and large-sized instances will be generated, the large-sized instances being of the size of the real instances of the OPS scheduling problem.

\begin{figure}[ht!]
\centering
\includegraphics[width=0.9\textwidth]{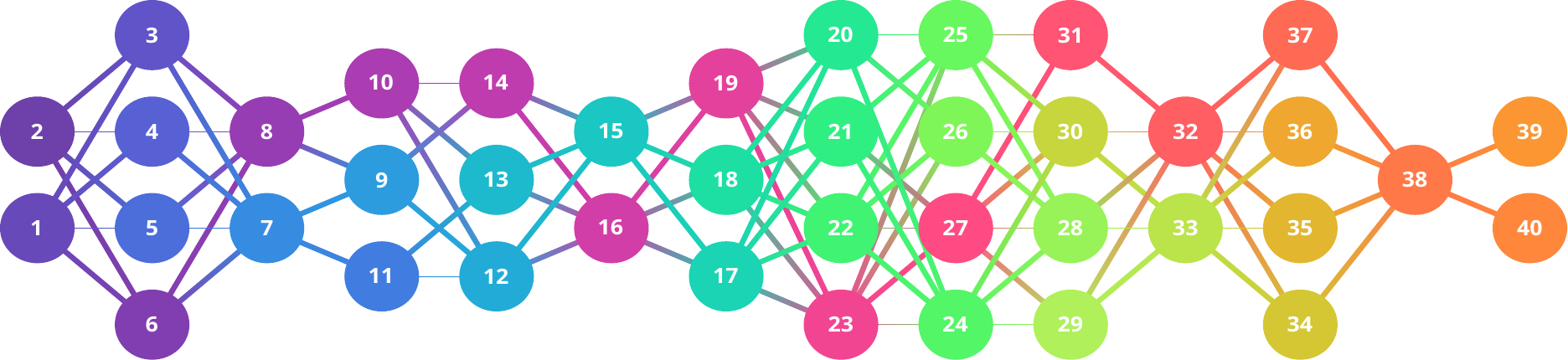} 
\caption{Random DAG representing the precedence constraints of a job with~40 operations. All arcs are directed from left to right. The visual representation of the DAG was drawn by using the network simplex layering proposed in~\cite{gansner1993technique}.}
\label{fig:random_ops_job}
\end{figure}

The generation of an instance relies on six given integer parameters, namely, the number of jobs~$n$, the minimum~$o_{\min}$ and maximum~$o_{\max}$ number of operations \textit{per} job, the minimum~$m_{\min}$ and the maximum~$m_{\max}$ number of machines, and the maximum number~$q$ of periods of unavailability \textit{per} machine. Precedence constraints are defined as follows; starting with~$A=\emptyset$. For each job $j \in \{1,\dots,n\}$, the generation of its associated DAG starts by choosing its number of operations~$o_j \in [o_{\min},o_{\max}]$; and it proceeds by layers. Starting by layer~$L_0$, between~$1$ and~$4$ operations are chosen to populate~$L_0$. Additional layers~$L_i$ are also populated with~$1$ to~$4$ operations, until the number of operations~$o_j$ is reached. Operations in~$L_0$ have no predecessors; while operations in the last layer have no successors. For every layer~$L_i$, unless the last one, and every~$v \in L_i$, an operation $w \in L_{i+1}$ is randomly chosen and the arc $(v,w)$ is added to~$A$. All the other pairs $(v,w)$ with $v \in L_i$ and $w \in L_{i+1}$ are included in~$A$ with probability~$0.85$. Figure~\ref{fig:random_ops_job} shows the random DAG of a job~$j$ with $o_j=40$. Note that the total number of operations is defined as $o = \sum_{j=1}^n o_j$. For each operation~$i$ such that an arc of the form $(i,\cdot) \in A$ exists, with probability~$0.1$, the overlapping coefficient~$\theta_i$ is a real number chosen in $[0.5,0.99]$, otherwise, $\theta_i$ is set to~$1$.

The number of machines is given by a random number $m \in [m_{\min},m_{\max}]$. For each operation~$i$, the cardinality of~$F(i)$ is given by a random number in $[\ceil{0.3m},\ceil{0.7m}]$ and the elements of $F(i) \subseteq \{1,\dots,m\}$ are randomly chosen. Then, a machine~$\hat k \in F(i)$ and the associated integer processing time $p_{i \hat k} \in [1,99]$ are randomly chosen. For all other machines $k \in F(i)$, $k \neq \hat k$, $p_{ik} \in [p_{i \hat k},\min\{3 p_{i \hat k},99\}]$ is randomly chosen. The number of periods of unavailability~$q_k$ of a machine~$k$ is chosen at random in~$[1,q]$. Let $\varphi_k$ be the mean of the processing times~$p_{ik}$ of the operations~$i$ such that $k \in F(i)$. Then, we define $a_k = 1 + \lceil \frac{\varphi_k}{q_k} \rceil$ as the distance between consecutive periods of unavailability. The first period of unavailability is given by $[\underline{u}_1^k,\bar u_1^k] = [a_k,\underline{u}_1^k + \lceil \frac{a_k}{R^k_1} \rceil + 1]$, where $R_1^k \in [2,10]$ is a random integer number. For $\ell=2,\dots,q_k$, the $\ell$-th period of unavailability is given by $[\underline{u}_{\ell}^k,\bar u_{\ell}^k] = [\bar u_{\ell-1}^k + a_k, \underline{u}_{\ell}^k + \lceil \frac{a_k}{R_{\ell}^k} \rceil + 1]$, where $R_{\ell}^k \in [2,10]$ is a random integer number.

Each operation~$i$ has three randomly chosen integer values $\mathit{size}_i \in [1,\bar S]$, $\mathit{color}_i \in [1,\bar C]$, and $\mathit{varnish}_i \in [1,\bar V]$, with $\bar S=10$, $\bar C=4$, and $\bar V=6$, that stand for the operation's size, color, and varnish, respectively. Consider two operations~$i$ and~$j$ that are processed consecutively, $i$ before~$j$, on a machine~$k$. If $\mathit{size}_i < \mathit{size}_j$, then $\mathit{st}'_k$ units of time are required to setup the machine; while $\mathit{st}''_k$ units of time are required if $\mathit{size}_i > \mathit{size}_j$. (No setup time is required, regarding the size feature, if $\mathit{size}_i = \mathit{size}_j$.) If $\mathit{color}_i \neq \mathit{color}_j$, $\mathit{ct}_k$ additional units of time are required; and if $\mathit{varnish}_i \neq \mathit{varnish}_j$, other additional $\mathit{vt}_k$ units of time are required. Values $\mathit{st}'_k$, $\mathit{st}''_k$, $\mathit{ct}_k$, and $\mathit{vt}_k$ are (machine dependent) random integer values in~$[2,6]$. The sum of the required values composes the setup time~$\gamma^I_{ijk} \geq 0$. The setup time $\gamma^F_{ik}$ for the case in which operation~$i$ is assigned to machine~$k$ and it is the first operation to be executed on the machine is given by $\gamma^F_{ik} = \max\{\mathit{st}'_k, \mathit{st}''_k\} + \mathit{ct}_k + \mathit{vt}_k$. It should be noted that the duration of a setup operation in between two consecutive operations~$i$ and~$j$ is not related to the processing time of the operations. A clear example of this situation corresponds to the setup operations of a cutting machine. The setup of the machine corresponds to adjusting the machine from the size of printed sheets and the cutting pattern of operation~$i$ to the size of printed sheets and the cutting pattern of operation~$j$; and this adjustment is not related to the quantity of pieces that must be cut in the two operations.

The release time~$r_i$ of each operation~$i$ is equal to~$0$ with probability~$0.975$. When~$r_i$ is not zero, it is a random integer number within the interval $[1,99]$. An operation~$i$ with no predecessors has probability~$0.01$ of belonging to~$T$; i.e.,\ having a fixed starting time~$\bar s_i$ at a predefined machine~$k \in F(i)$. (At most one fixed operation per machine is allowed.) If this is the case, $p_{ik}$ is redefined as a random number in~$[1,99]$, $F(i)$ is redefined as the singleton $F(i)=\{k\}$, and $\bar s_i$ is randomly chosen in~$[\gamma^F_{\mathrm{ub}}, \underline{u}_1^k - p_{ik}]$, where $\gamma^F_{\mathrm{ub}}= \bar S + \bar C + \bar V$ is an upper bound of the setup time of an operation that is the first operation to be processed by a machine, as it is the case of a fixed operation generated in this way. 

\subsection{Numerical results} 
\label{secres}

The experiments, carried out using the High-Performance Computing (HPC) facilities of the University of Luxembourg~\citep{VBCG_HPCS14}, were conducted on an Intel Xeon E5-2680 v4 2.4 GHz with 4GB memory (per core) running CentOS Linux~7.7 (in 64-bit mode); always using a single physical core.

By default, see~\cite[pp.~227--228]{cplexmanual}, a solution to a MILP model is reported by CPLEX when \[
\text{absolute gap} = \text{incumbent solution} - \text{best lower bound} \leq \epsilon_\text{abs}
\]
or
\[
\text{relative gap} = \frac{|\text{incumbent solution} - \text{best lower bound}|}{10^{-10} + |\text{incumbent solution}|} \leq \epsilon_\text{rel},
\]
with $\epsilon_\text{abs} = 10^{-6}$ and $\epsilon_\text{rel} = 10^{-4}$. In the considered instances, the objective function~(\ref{eq:min_cmax}) assumes integer values at feasible points. Thus, on the one hand, a stopping criterion based on a relative error less than or equal to $\epsilon_\text{rel} = 10^{-4}$ may have the undesired effect of stopping the method prematurely; and, on the other hand, an absolute error strictly smaller than 1 is enough to prove the optimality of the incumbent solution. Therefore, following~\cite{bromro} and~\cite{abmro2}, we considered $\epsilon_\text{abs} = 1 - 10^{-6}$ and $\epsilon_\text{rel} = 0$. All other parameters of CPLEX were kept with their default values. CP Optimizer was run with all its default parameters. A CPU time limit of two hours per instance was imposed.

\subsubsection{Experiments with small-sized instances}

\begin{figure}[ht!]
\centering
\includegraphics[width=\linewidth]{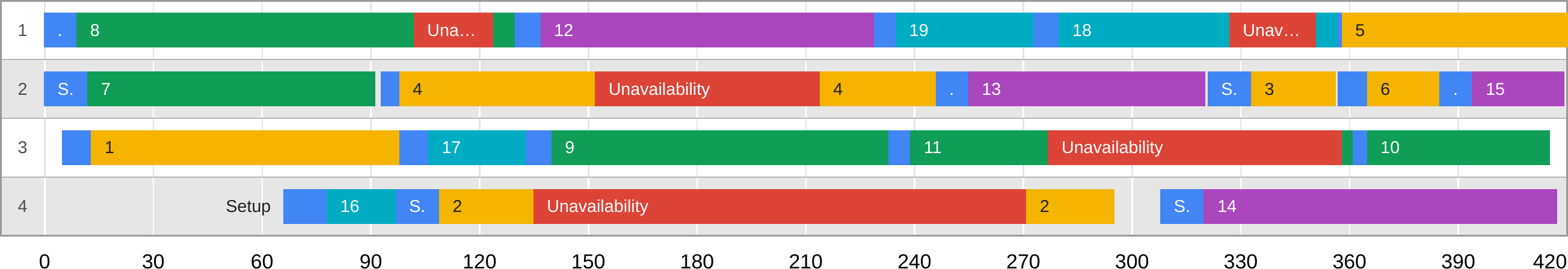}
\caption{Optimal solution to instance~30 of the set of small-sized instances. Operations with the same color belong to the same job; while setups are represented in blue and machines' unavailiability periods in red. The complete instance's data can be found in \url{https://willtl.github.io/ops}.}
\label{solsmall}
\end{figure}

In this section, a set of~30 small-sized instances is considered. The $k$-th instance was generated with the following parameters: $n = 1 + \lceil \frac{k}{30} \times 3 \rceil$, $o_{\min} = 2$, $o_{\max} = 3 + \lceil \frac{k}{30} \times 2 \rceil$, $m_{\min} = 2$, $m_{\max} = 3 + \lceil \frac{k}{30} \times 2 \rceil$, and $q = 4$. Table~\ref{tab1} shows the main characteristics of each instance; while Table~\ref{tab2} shows the value of the solutions found when solving the MILP and the CP Optimizer formulations with IBM ILOG CPLEX and IBM ILOG CP Optimizer solvers, respectively. The associated effort measurements are also shown in the table. Most of the columns in the table are self-explanatory. When optimality is not proven, column ``Makespan'' shows the best lower bound, the best upper bound, and the gap. The precise CP Optimizer model that is being solved and the meaning of the columns ``number of branches in phases~1 and~2'' will be elucidated in the next section. The CPU time is expressed in seconds. Figures in the table clearly show that the CP Optimizer solver outperformed the MILP solver. In the ten instances with $n=2$, both solvers performed similarly. In the ten instances with $n=3$, the CP Optimizer solver outperformed the MILP solver by one or two orders of magnitude. In the ten instances with $n=4$, the MILP solver was not able to prove optimality of any of the instances; while the CP Optimizer solver solved all instances in a few seconds of CPU time. As an illustration, Figure~\ref{solsmall} shows the solution to instance~30 found for the MILP and the CP Optimizer formulations by IBM ILOG CPLEX and IBM ILOG CP Optimizer solvers, respectively, but proven to be optimal in the latter case only.

\begin{table}[ht!]
\caption{Main features of the considered thirty small-sized instances.}
\label{tab1}
\centering 
\resizebox{\textwidth}{!}{
\begin{tabular}{cccccccccccc}
\toprule
& \multicolumn{6}{c}{Main instance characteristics} 
& \multicolumn{3}{c}{MILP formulation} & \multicolumn{2}{c}{CP Optimizer formulation} \\
\midrule
\multirow{2}{*}{Instance} & 
\multirow{2}{*}{$m$} & 
\multirow{2}{*}{$\sum_{k=1}^m q_k$} & 
\multirow{2}{*}{$n$} & 
\multirow{2}{*}{$o$} & 
\multirow{2}{*}{$|A|$} & 
\multirow{2}{*}{$|T|$} & 
\# binary & \# continuous & \multirow{2}{*}{\# constraints} & \# integer & \multirow{2}{*}{\# constraints} \\
& & & & & & & variables & variables & & variables \\
\midrule 
1 & 3 & 7 & 2 & 9 & 10 & 1 & 248 & 109 & 889 & 76 & 209 \\ \rowcolor{gray!10}
2 & 4 & 10 & 2 & 8 & 8 & 1 & 220 & 99 & 784 & 67 & 187 \\ 
3 & 3 & 8 & 2 & 8 & 7 & 1 & 238 & 99 & 843 & 69 & 191 \\ \rowcolor{gray!10}
4 & 4 & 10 & 2 & 9 & 9 & 1 & 263 & 113 & 928 & 77 & 217 \\ 
5 & 2 & 5 & 2 & 9 & 8 & 2 & 200 & 101 & 752 & 58 & 162 \\ \rowcolor{gray!10}
6 & 2 & 4 & 2 & 8 & 7 & 1 & 184 & 93 & 685 & 57 & 159 \\ 
7 & 2 & 6 & 2 & 8 & 6 & 0 & 234 & 93 & 832 & 54 & 154 \\ \rowcolor{gray!10}
8 & 4 & 11 & 2 & 7 & 6 & 0 & 240 & 93 & 836 & 69 & 192 \\ 
9 & 4 & 8 & 2 & 8 & 7 & 1 & 212 & 103 & 766 & 81 & 214 \\ \rowcolor{gray!10}
10 & 4 & 13 & 2 & 10 & 8 & 0 & 425 & 131 & 1416 & 104 & 280 \\ 
11 & 4 & 8 & 3 & 14 & 14 & 0 & 549 & 181 & 1805 & 135 & 372 \\ \rowcolor{gray!10}
12 & 3 & 9 & 3 & 11 & 8 & 0 & 354 & 131 & 1225 & 83 & 231 \\ 
13 & 3 & 7 & 3 & 14 & 13 & 0 & 458 & 163 & 1546 & 100 & 281 \\ \rowcolor{gray!10}
14 & 3 & 8 & 3 & 13 & 12 & 0 & 514 & 163 & 1701 & 116 & 323 \\ 
15 & 2 & 4 & 3 & 13 & 12 & 1 & 322 & 145 & 1169 & 83 & 236 \\ \rowcolor{gray!10}
16 & 4 & 8 & 3 & 13 & 12 & 0 & 485 & 169 & 1619 & 121 & 341 \\ 
17 & 2 & 4 & 3 & 14 & 14 & 1 & 301 & 151 & 1131 & 77 & 224 \\ \rowcolor{gray!10}
18 & 4 & 7 & 3 & 12 & 11 & 2 & 339 & 151 & 1182 & 109 & 300 \\ 
19 & 4 & 13 & 3 & 14 & 14 & 1 & 669 & 183 & 2165 & 134 & 380 \\ \rowcolor{gray!10}
20 & 4 & 11 & 3 & 12 & 9 & 1 & 466 & 153 & 1559 & 107 & 300 \\ 
21 & 4 & 10 & 4 & 18 & 16 & 1 & 927 & 237 & 2936 & 172 & 490 \\ \rowcolor{gray!10}
22 & 2 & 3 & 4 & 20 & 18 & 0 & 628 & 223 & 2149 & 126 & 359 \\ 
23 & 3 & 8 & 4 & 18 & 15 & 0 & 895 & 225 & 2824 & 144 & 419 \\ \rowcolor{gray!10}
24 & 3 & 6 & 4 & 18 & 14 & 0 & 822 & 227 & 2596 & 159 & 447 \\ 
25 & 2 & 5 & 4 & 18 & 17 & 2 & 615 & 201 & 2097 & 114 & 326 \\ \rowcolor{gray!10}
26 & 2 & 6 & 4 & 19 & 21 & 2 & 622 & 207 & 2159 & 112 & 321 \\ 
27 & 2 & 2 & 4 & 23 & 26 & 1 & 624 & 249 & 2217 & 127 & 374 \\ \rowcolor{gray!10}                
28 & 3 & 7 & 4 & 19 & 17 & 1 & 1027 & 241 & 3192 & 157 & 459 \\             
29 & 3 & 10 & 4 & 18 & 18 & 1 & 1057 & 231 & 3282 & 153 & 449 \\ \rowcolor{gray!10}
30 & 4 & 12 & 4 & 19 & 19 & 1 & 1132 & 255 & 3516 & 183 & 530 \\ 
\bottomrule
\end{tabular}}  
\end{table}

\begin{table}[ht!]
\centering 
\caption{Description of the solutions found and effort measurements of the IBM ILOG CPLEX and IBM ILOG CP Optimizer applied to the thirty small-sized instances.}
\label{tab2}
\resizebox{0.9\textwidth}{!}{\begin{tabular}{ccccccccc}
\toprule
\multirow{4}{*}{}  & \multicolumn{4}{c}{IBM ILOG CPLEX}& \multicolumn{4}{c}{IBM ILOG CP Optimizer} \\
\cmidrule(lr){2-9}
& \multirow{3}{*}{Makespan} & \multicolumn{3}{c}{Effort measurement} & \multirow{3}{*}{Makespan} & \multicolumn{3}{c}{Effort measurement} \\
\cmidrule(lr){3-5} \cmidrule{7-9}
& & \multirow{2}{*}{\# iterations} & \multirow{2}{*}{\# B\&B Nodes} & \multirow{2}{*}{CPU} & & \multicolumn{2}{c}{\# of branches in} & \multirow{2}{*}{CPU} \\
& & & & & & phase 1 & phase 2 & \\
\midrule
1 & 274  & 9286 & 824 & 0.5 & 274 & 299 & 21 & 0.1 \\ \rowcolor{gray!10}
2 & 230  & 7302 & 705 & 0.4 & 230 & 215 & 3 & 0.1 \\
3 & 337  & 2451 & 231 & 0.2 & 337 & 66 & 3 & 0.1 \\ \rowcolor{gray!10}
4 & 276  & 6410 & 444 & 0.4 & 276 & 179 & 21 & 0.1 \\
5 & 495  & 6444 & 763 & 0.3 & 495 & 212 & 21 & 0.1 \\ \rowcolor{gray!10}
6 & 271  & 5645 & 665 & 0.3 & 271 & 193 & 62 & 0.1 \\
7 & 370  & 7223 & 561 & 0.3 & 370 & 203 & 580 & 0.1 \\ \rowcolor{gray!10}
8 & 279  & 3263 & 114 & 0.3 & 279 & 98 & 17 & 0.1 \\
9 & 274  & 1261 & 113 & 0.1 & 274 & 8 & 19 & 0.1 \\ \rowcolor{gray!10}
10 & 329  & 8530 & 357 & 0.7 & 329 & 1148 & 30 & 0.1 \\
11 & 239  & 649022 & 21705 & 56.6 & 239 & 5274 & 3 & 0.2 \\ \rowcolor{gray!10}
12 & 273  & 54160 & 2870 & 3.9 & 273 & 410 & 3 & 0.1 \\
13 & 266  & 2464645 & 79179 & 176.0 & 266 & 14835 & 10038 & 1.0 \\ \rowcolor{gray!10}
14 & 518  & 401863 & 20863 & 42.4 & 518 & 1669 & 29 & 0.1 \\
15 & 551  & 1826080 & 75865 & 153.8 & 551 & 14824 & 95573 & 1.2 \\ \rowcolor{gray!10}
16 & 278  & 21274 & 1340 & 1.4 & 278 & 361 & 29 & 0.1 \\
17 & 540  & 108667 & 8321 & 7.7 & 540 & 14516 & 3 & 0.2 \\ \rowcolor{gray!10}
18 & 327  & 4623 & 441 & 0.3 & 327 & 12 & 27 & 0.1 \\
19 & 325  & 248850 & 11251 & 16.5 & 325 & 1665 & 990 & 0.1 \\ \rowcolor{gray!10}
20 & 264  & 16287 & 826 & 1.1 & 264 & 1500 & 179 & 0.1 \\
21 & [263, 300] 12.3\% & 50935466 & 1465096 & 7200 & 300 & 21643 & 11865 & 0.7 \\ \rowcolor{gray!10}
22 & [204, 671] 69.6\% & 73458018 & 1400247 & 7200 & 651 & 220130 & 3 & 8.0 \\
23 & [406, 467] 13.1\% & 56974881 & 1217720 & 7200 & 467 & 19771 & 3 & 0.6 \\ \rowcolor{gray!10}
24 & [336, 572] 41.3\% & 80022693 & 1673003 & 7200 & 571 & 72971 & 3 & 2.6 \\
25 & [436, 672] 35.1\% & 62426552 & 2011692 & 7200 & 672 & 76887 & 110821 & 3.0 \\ \rowcolor{gray!10}
26 & [374, 628] 40.4\% & 88226422 & 1926203 & 7200 & 627 & 226752 & 324106 & 12.5 \\
27 & [223, 719] 69.0\% & 79180995 & 1529405 & 7200 & 702 & 498600 & 927317 & 63.7 \\ \rowcolor{gray!10}
28 & [297, 493] 39.8\% & 62733982 & 707445 & 7200 & 437 & 28733 & 3 & 1.4 \\
29 & [343, 498] 31.1\% & 55996388 & 1397887 & 7200 & 480 & 100627 & 39 & 4.0 \\ \rowcolor{gray!10}
30 & [417, 420] 0.7\% & 65757996 & 1158447 & 7200 & 420 & 8148 & 10885 & 0.3 \\
\midrule
Mean & 405.2 \; 11.75\% & 22718889.3 & 490486.1 & 2415.44 & 401.43 \; 0.0\% & 44398.3 & 49756.53 & 3.37 \\
\bottomrule
\end{tabular}}
\end{table} 

\subsubsection{Experiments with medium-sized instances} 

In this section, a set of~20 medium-sized instances is considered. The $k$-th instance was generated with the following parameters: $n = 4 + \lceil \frac{k}{20} \times 6 \rceil$, $o_{\min} = 6$, $o_{\max} = 7 + \lceil \frac{k}{20} \times 5 \rceil$, $m_{\min} = 6$, $m_{\max} = 7 + \lceil \frac{k}{20} \times 13 \rceil$, and $q = 8$. Table~\ref{tab3} shows the main characteristics of each instance.

\begin{figure}[ht!]
\centering
\includegraphics[width=\linewidth]{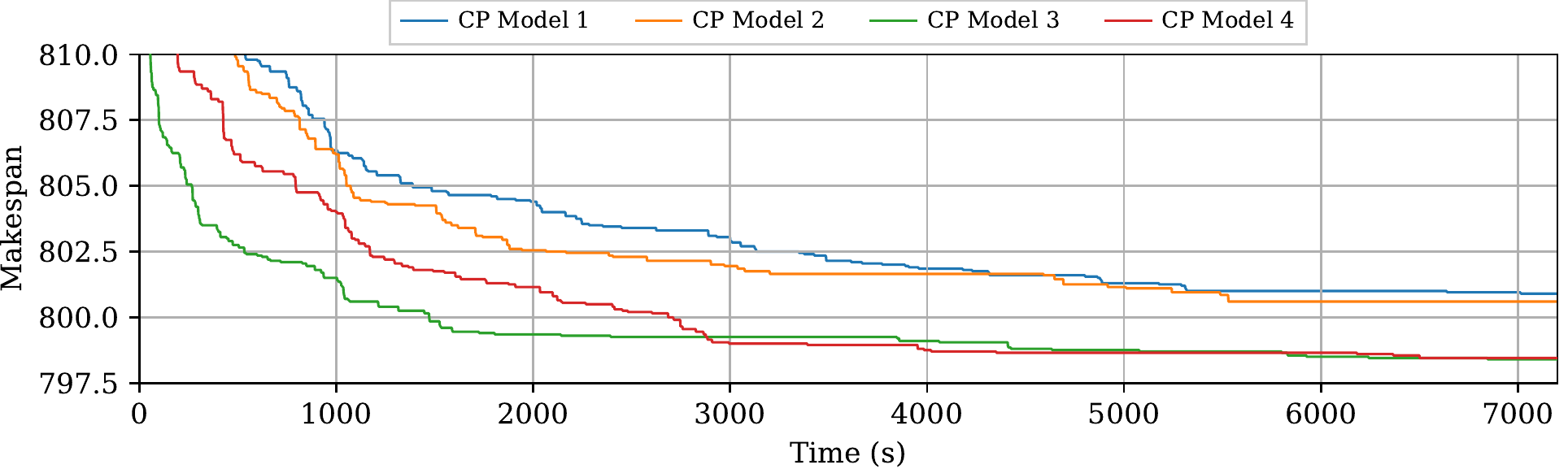}
\caption{Evolution over time of the (average of the) incumbent solutions' makespan of medium-sized instances along the resolution of CP Models~1--4.}
\label{evolution}
\end{figure}

As described at the end of Section~\ref{sec:cp}, the CP Optimizer formulation of the OPS scheduling problem is given by the minimization of~(\ref{eqcpo:3}) subject to constraints~(\ref{eqcpo:5}--\ref{eqcpo:1},\ref{eqcpoops:B2},\ref{eqcpoops:1}--\ref{eqcpoops:6bis},\ref{eqcpoops:10}--\ref{eqcpoops:20}); with  constraints~(\ref{eqcpo:7}) and~(\ref{eqcpoops:12b},\ref{eqcpoops:18}) being optional. Therefore, in a first experiment, we aim to evaluate the influence of the optional constraints by comparing: (i) the plain model (named CP Model~1 from on); (ii) the plain model plus constraints~(\ref{eqcpo:7}) (named CP Model~2 from on); (iii) the plain model plus constraints~(\ref{eqcpoops:12b},\ref{eqcpoops:18}) (named CP Model~3 from on); and (iv) the plain model plus constraints~(\ref{eqcpo:7}) and constraints~(\ref{eqcpoops:12b},\ref{eqcpoops:18}) (named CP Model~4 from on). Table~\ref{tab4} shows the value of the solutions found when solving each of the four CP Optimizer models with IBM ILOG CP Optimizer solver. The associated effort measurements are also shown in the table. Figures in the table show that considering the optional constraints~(\ref{eqcpoops:12b},\ref{eqcpoops:18}) helps CP Optimizer solver to close the gap and prove optimality in 12 out of the~20 considered instances when solving CP Models~3 and~4; while, when solving CP Models~1 and~2, that do not consider the optional constraints~(\ref{eqcpoops:12b},\ref{eqcpoops:18}), gaps are closed in~10 out of the~20 considered instances. On the other hand, including or not the optional constraints~(\ref{eqcpo:7}) appears to have no relevant influence on the resolution process of the considered set of instances. Figure~\ref{evolution} shows the average value of the makespan, over the~20 considered instances, as a function of time, over the resolution of CP Models~1--4. The graphic confirms, as expected, that, CP Optimizer solver is able to improve the incumbent solution faster when applied to CP Models~3 and~4 compared to its application to CP Models~1 and~2.

\begin{table}[ht!]
\centering 
\caption{Main features of the considered twenty medium-sized instances.}
\label{tab3}
\resizebox{\textwidth}{!}{\begin{tabular}{cccccccccccc}
\toprule
& \multicolumn{6}{c}{Main instance characteristics} 
& \multicolumn{3}{c}{MILP formulation} & \multicolumn{2}{c}{CP Optimizer formulation} \\
\midrule
\multirow{2}{*}{Instance} & 
\multirow{2}{*}{$m$} & 
\multirow{2}{*}{$\sum_{k=1}^m q_k$} & 
\multirow{2}{*}{$n$} & 
\multirow{2}{*}{$o$} & 
\multirow{2}{*}{$|A|$} & 
\multirow{2}{*}{$|T|$} & 
\# binary & \# continuous & \multirow{2}{*}{\# constraints} & \# integer & \multirow{2}{*}{\# constraints} \\
& & & & & & & variables & variables & & variables \\
\midrule 
1 & 8 & 41 & 5 & 39 & 58 & 0 & 5987 & 635 & 16984 & 584 & 1671 \\  \rowcolor{gray!10}
2 & 7 & 28 & 5 & 36 & 54 & 3 & 4341 & 559 & 12398 & 471 & 1381 \\ 
3 & 8 & 46 & 5 & 43 & 75 & 1 & 7753 & 717 & 21670 & 643 & 1895 \\  \rowcolor{gray!10}
4 & 7 & 27 & 6 & 43 & 60 & 1 & 6615 & 693 & 18246 & 595 & 1749 \\ 
5 & 6 & 31 & 6 & 45 & 61 & 3 & 6189 & 659 & 17641 & 512 & 1519 \\  \rowcolor{gray!10}
6 & 7 & 38 & 6 & 46 & 60 & 0 & 7292 & 715 & 20646 & 579 & 1725 \\ 
7 & 8 & 31 & 7 & 64 & 108 & 1 & 12302 & 1049 & 33289 & 929 & 2731 \\  \rowcolor{gray!10}
8 & 9 & 51 & 7 & 53 & 76 & 1 & 11814 & 933 & 32185 & 860 & 2524 \\ 
9 & 6 & 25 & 7 & 56 & 90 & 2 & 9139 & 839 & 25240 & 677 & 2011 \\  \rowcolor{gray!10}
10 & 8 & 31 & 7 & 63 & 110 & 2 & 12993 & 1057 & 34805 & 932 & 2777 \\ 
11 & 13 & 52 & 8 & 75 & 118 & 0 & 26493 & 1615 & 67250 & 1693 & 4962 \\  \rowcolor{gray!10}
12 & 11 & 56 & 8 & 78 & 142 & 2 & 26370 & 1535 & 68106 & 1554 & 4543 \\ 
13 & 16 & 75 & 8 & 68 & 105 & 0 & 30571 & 1713 & 77029 & 1949 & 5652 \\  \rowcolor{gray!10}
14 & 12 & 61 & 9 & 72 & 105 & 0 & 26915 & 1533 & 68895 & 1610 & 4685 \\ 
15 & 17 & 73 & 9 & 76 & 107 & 0 & 37144 & 1965 & 92371 & 2239 & 6508 \\  \rowcolor{gray!10}
16 & 15 & 67 & 9 & 89 & 156 & 2 & 42871 & 2113 & 106337 & 2287 & 6756 \\ 
17 & 7 & 34 & 10 & 109 & 207 & 1 & 32729 & 1723 & 86417 & 1460 & 4362 \\  \rowcolor{gray!10}
18 & 17 & 72 & 10 & 96 & 164 & 1 & 53876 & 2457 & 131704 & 2827 & 8190 \\ 
19 & 14 & 72 & 10 & 92 & 142 & 2 & 44376 & 2095 & 110957 & 2273 & 6617 \\  \rowcolor{gray!10}
20 & 11 & 51 & 10 & 91 & 135 & 0 & 36772 & 1843 & 92819 & 1850 & 5439 \\
\bottomrule
\end{tabular}}
\end{table}

\begin{table}[ht!]
\centering 
\caption{Description of the solutions found and effort measurements of the IBM ILOG CP Optimizer applied to the twenty medium-sized instances.}
\label{tab4}
\resizebox{\textwidth}{!}{\begin{tabular}{ccccccccccccc}
\toprule
\multirow{3}{*}{} & 
\multicolumn{3}{c}{CP Model 1} & \multicolumn{3}{c}{CP Model 2} & 
\multicolumn{3}{c}{CP Model 3} & \multicolumn{3}{c}{CP Model 4} \\
\cmidrule{2-13}
& 
\multirow{2}{*}{Makespan} & \multicolumn{2}{c}{Effort measurement} & 
\multirow{2}{*}{Makespan} & \multicolumn{2}{c}{Effort measurement} & 
\multirow{2}{*}{Makespan} & \multicolumn{2}{c}{Effort measurement} & 
\multirow{2}{*}{Makespan} & \multicolumn{2}{c}{Effort measurement} \\
\cmidrule{3-4} \cmidrule{6-7} \cmidrule{9-10} \cmidrule{12-13} 
&
& \# branches & CPU & & \# branches & CPU & 
& \# branches & CPU & & \# branches & CPU \\
\midrule 
1 & 344 & 11819616 & 599.7 & 344 & 4945664 & 499.0 & 344 & 368683 & 16.1 & 344 & 399038 & 32.5 \\ \rowcolor{gray!10}
2 & 357 & 129666 & 12.0 & 357 & 112211 & 11.3 & 357 & 35850 & 3.9 & 357 & 34779 & 3.8 \\
3 & 404 & 45163467 & 4315.8 & [361, 409] 11.7\% & 77088019 & 7200 & [361, 406] 11.1\% & 88706470 & 7200 & 404 & 18733373 & 2825.1 \\ \rowcolor{gray!10}
4 & 458 & 17065512 & 1081.0 & 458 & 7777910 & 531.7 & 458 & 269785 & 25.4 & 458 & 548398 & 49.5 \\
5 & [474, 510] 7.1\% & 112252005 & 7200 & [474, 515] 8.0\% & 86834261 & 7200 & 506 & 384590 & 40.8 & 506 & 399229 & 90.5 \\ \rowcolor{gray!10}
6 & [329, 454] 27.5\% & 86652702 & 7200 & [334, 438] 23.7\% & 84114310 & 7200 & [334, 437] 23.6\% & 80406442 & 7200 & [335, 442] 24.2\% & 47105274 & 7200 \\
7 & 2429 & 1262 & 0.2 & 2429 & 1272 & 0.2 & 2429 & 1455 & 0.2 & 2429 & 1488 & 0.2 \\ \rowcolor{gray!10}
8 & [360, 456] 21.1\% & 89059560 & 7200 & [360, 451] 20.2\% & 64334715 & 7200 & [360, 460] 21.7\% & 65621010 & 7200 & [360, 452] 20.4\% & 59032373 & 7200 \\
9 & [629, 631] 0.3\% & 118721872 & 7200 & [629, 636] 1.1\% & 34808015 & 7200 & [629, 630] 0.2\% & 94392736 & 7200 & [629, 630] 0.2\% & 86015687 & 7200 \\ \rowcolor{gray!10}
10 & 1184 & 1079 & 0.1 & 1184 & 1052 & 0.1 & 1184 & 1204 & 0.2 & 1184 & 1189 & 0.3 \\
11 & [406, 430] 5.6\% & 39491824 & 7200 & [406, 431] 5.8\% & 53721744 & 7200 & [406, 427] 4.9\% & 48553125 & 7200 & [406, 424] 4.2\% & 45050481 & 7200 \\ \rowcolor{gray!10}
12 & [457, 508] 10.0\% & 75583273 & 7200 & [457, 510] 10.4\% & 76775244 & 7200 & [457, 503] 9.1\% & 57651439 & 7200 & [457, 505] 9.5\% & 48316050 & 7200 \\
13 & 347 & 148874 & 26.5 & 347 & 144200 & 15.8 & 347 & 64381 & 7.7 & 347 & 29387 & 6.6 \\ \rowcolor{gray!10}
14 & [302, 412] 26.7\% & 43670920 & 7200 & [320, 404] 20.8\% & 57130518 & 7200 & [302, 399] 24.3\% & 48125327 & 7200 & [320, 403] 20.6\% & 33552405 & 7200 \\
15 & 319 & 259273 & 40.6 & 319 & 275642 & 58.8 & 319 & 55919 & 12.4 & 319 & 94552 & 27.6 \\ \rowcolor{gray!10}
16 & 543 & 138290 & 18.0 & 543 & 1814 & 0.5 & 543 & 1343 & 0.4 & 543 & 1331 & 0.8 \\
17 & [1044, 1053] 0.9\% & 60769510 & 7200 & [1044, 1059] 1.4\% & 48747654 & 7200 & 1052 & 12862072 & 1589.7 & [1052, 1055] 0.3\% & 32700611 & 7200 \\ \rowcolor{gray!10}
18 & 3184 & 1968 & 0.6 & 3184 & 1968 & 0.6 & 3184 & 2091 & 0.8 & 3184 & 2091 & 0.7 \\
19 & [1449, 1451] 0.1\% & 93263868 & 7200 & [1449, 1451] 0.1\% & 83339741 & 7200 & 1451 & 2305 & 0.7 & 1451 & 2371 & 0.7 \\ \rowcolor{gray!10}
20 & [360, 544] 33.8\% & 71372646 & 7200 & [417, 543] 23.2\% & 51052068 & 7200 & [360, 532] 32.3\% & 43254587 & 7200 & [417, 532] 21.6\% & 41945390 & 7200 \\
\midrule
Mean & 800.9 \; 6.65\% & 43278359 & 3904.72 & 800.6 \; 6.32\% & 36560401 & 4015.9 & 798.4 \; 6.36\% & 27038040 & 2964.91 & 798.45 \; 5.05\% & 20698274 & 3031.91 \\ 
\bottomrule
\end{tabular}}
\end{table}  

\begin{table}[ht!]
\centering 
\caption{Description of the solutions found and effort measurements of the IBM ILOG CP Optimizer applied to the twenty medium-sized instances. The left half of the table corresponds to CP Model~4 (already shown in Table~\ref{tab4}); while the right half of the table corresponds to solving first the Incomplete CP Model and passing its solution as an initial guess for the resolution of CP Model~4.}
\label{tab5}
\resizebox{\textwidth}{!}{\begin{tabular}{cccccccccc}
\toprule
\multirow{4}{*}{} & 
\multicolumn{3}{c}{CP Model 4} & 
\multicolumn{6}{c}{Incomplete CP Model + CP Model 4} \\
\cmidrule{2-10}
& 
\multirow{4}{*}{Makespan} & \multicolumn{2}{c}{Effort measurement} & 
\multirow{4}{*}{Makespan} & \multicolumn{3}{c}{Effort measurement} \\
\cmidrule{3-4} \cmidrule{6-10}
& & \multirow{2}{*}{\# branches} & \multirow{2}{*}{CPU} & & \multicolumn{2}{c}{phase 1} & \multicolumn{2}{c}{phase 2} & Total \\
\cmidrule{6-7} \cmidrule{8-9}
& & & & & \# branches & CPU & \# branches & CPU & CPU \\
\midrule 
1 & 344 & 399038 & 32.5 & 344 & 288898 & 15.53 & 280015 & 15.07 & 30.6 \\ \rowcolor{gray!10}
2 & 357 & 34779 & 3.8 & 357 & 22228 & 1.09 & 3 & 0.01 & 1.1 \\
3 & 404 & 18733373 & 2825.1 & 404 & 132017215 & 4969.69 & 16246793 & 611.61 & 5581.3 \\ \rowcolor{gray!10}
4 & 458 & 548398 & 49.5 & 458 & 146895 & 5.59 & 3 & 0.01 & 5.6 \\
5 & 506 & 399229 & 90.5 & 506 & 272817 & 31.26 & 294203 & 33.74 & 65.0 \\ \rowcolor{gray!10}
6 & [335, 442] 24.2\% & 47105274 & 7200 & [335, 441] 24.0\% & 142959186 & 4800.0 & 21538650 & 2400.0 & 7200 \\
7 & 2429 & 1488 & 0.2 & 2429 & 1241 & 0.08 & 131 & 0.02 & 0.1 \\ \rowcolor{gray!10}
8 & [360, 452] 20.4\% & 59032373 & 7200 & [360, 450] 20.0\% & 129153428 & 4800.0 & 19362459 & 2400.0 & 7200 \\
9 & [629, 630] 0.2\% & 86015687 & 7200 & 629 & 43101316 & 1096.49 & 115 & 0.01 & 1096.5 \\ \rowcolor{gray!10}
10 & 1184 & 1189 & 0.3 & 1184 & 1250 & 0.09 & 3 & 0.01 & 0.1 \\
11 & [406, 424] 4.2\% & 45050481 & 7200 & [406, 418] 2.9\% & 123558867 & 4800.0 & 18454514 & 2400.0 & 7200 \\ \rowcolor{gray!10}
12 & [457, 505] 9.5\% & 48316050 & 7200 & [457, 499] 8.4\% & 128774356 & 4800.0 & 17042168 & 2400.0 & 7200 \\
13 & 347 & 29387 & 6.6 & 347 & 6914 & 0.39 & 3 & 0.01 & 0.4 \\ \rowcolor{gray!10}
14 & [320, 403] 20.6\% & 33552405 & 7200 & [320, 394] 18.8\% & 109752381 & 4800.0 & 13032611 & 2400.0 & 7200 \\
15 & 319 & 94552 & 27.6 & 319 & 47107 & 3.19 & 3 & 0.01 & 3.2 \\ \rowcolor{gray!10}
16 & 543 & 1331 & 0.8 & 543 & 1481 & 0.29 & 3 & 0.01 & 0.3 \\
17 & [1052, 1055] 0.3\% & 32700611 & 7200 & 1052 & 572991 & 59.72 & 455352 & 47.48 & 107.2 \\ \rowcolor{gray!10}
18 & 3184 & 2091 & 0.7 & 3184 & 2063 & 0.39 & 3 & 0.01 & 0.4 \\
19 & 1451 & 2371 & 0.7 & 1451 & 2504 & 0.36 & 187 & 0.04 & 0.4 \\ \rowcolor{gray!10}
20 & [417, 532] 21.6\% & 41945390 & 7200 & [417, 520] 19.8\% & 106764402 & 4800.0 & 12752639 & 2400.0 & 7200 \\ 
\midrule
Mean & 798.45 \; 5.05\% & 20698274 & 3031.91 & 796.45 \; 4.7\% & 45872377 & 1749.21 & 5972992 & 755.4 & 2504.61 \\ 
\bottomrule
\end{tabular}}
\end{table}  

\begin{figure}[ht!]
\centering
\includegraphics[width=\linewidth]{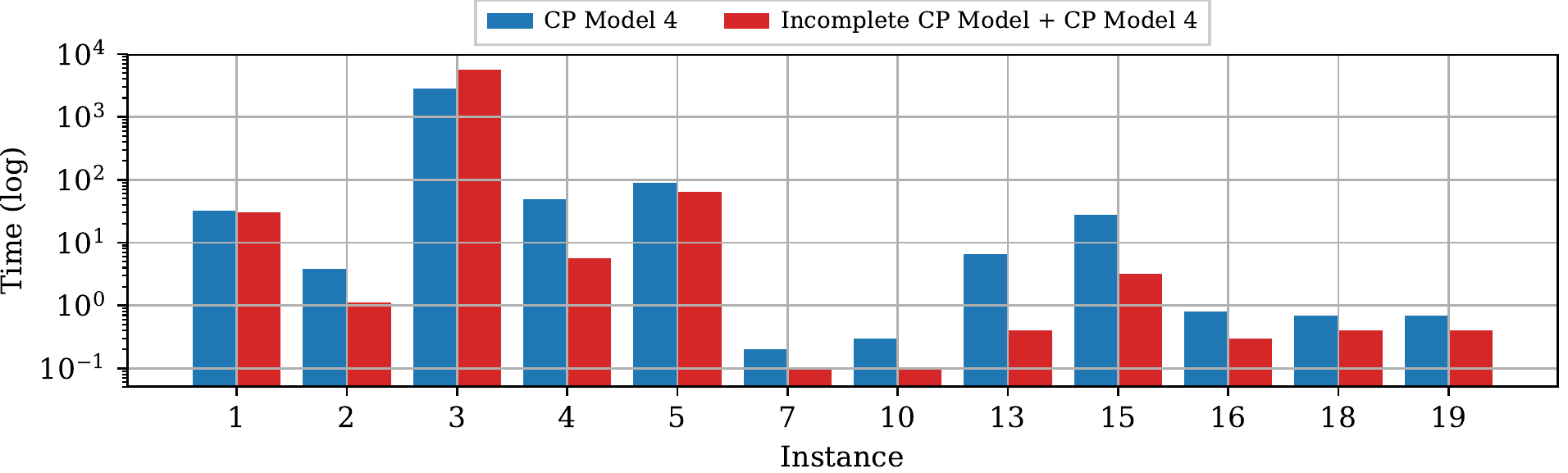}
\caption{Time comparison between the resolution of the CP Model~4 and its resolution in two phases considering the solution found to the Incomplete CP Model as an initial guess.}
\label{fig:2cpmodels}
\end{figure} 

\begin{figure}
\centering
\includegraphics[width=\linewidth]{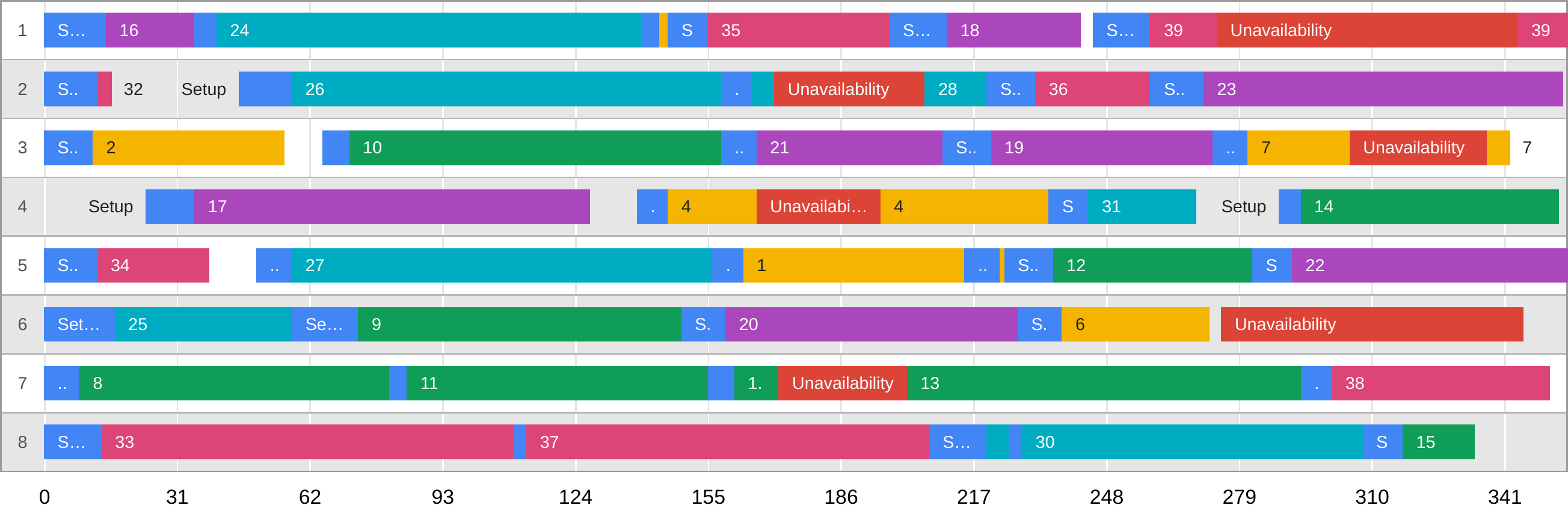}
\caption{Optimal solution to instance~1 of the set of medium-sized instances. Operations with the same color belong to the same job; while setups are represented in blue and machines' unavailiability periods in red. The complete instance's data can be found in \url{https://willtl.github.io/ops}.}
\label{fig:mops1}
\end{figure}

As described in Section~\ref{sec:setup}, modeling a setup in between consecutive operations is relatively easy in the CP Optimizer language. However, the simple formulation given by constraints~(\ref{eqcpoops:12b},\ref{eqcpoops:18}) does not consider the setup time of the first operation processed by each machine; as well as it does not consider that setup operations are not resumable, i.e.,\ they can not be interrupted by periods of unavailability of the machines. This is why the model that consists in minimizing~(\ref{eqcpo:3}) subject to constraints~(\ref{eqcpo:5}--\ref{eqcpo:1},\ref{eqcpoops:B2},\ref{eqcpoops:1}--\ref{eqcpoops:18}) can be considered an incomplete model of the OPS scheduling problem, and we name it Incomplete CP Model from now on. The Incomplete CP Model is much simpler that CP Model~4 and, as a resolution strategy, it can be solved first; and its optimal (or best known feasible) solution used as an initial guess to the resolution of CP Model~4. This two-phases strategy is based on the fact that IBM ILOG CP Optimizer solver has the capability of accepting as initial guess a possible infeasible and incomplete solution. By incomplete we mean that CP Model~4 has variables that are not present in the Incomplete CP Model; and, by infeasible, we mean that, with very high probability, a solution to the Incomplete CP Model has setup operations being interrupted by machines' unavailabilities, as well as it does not consider the setup operation of the first operation being processed by each machine; thus being infeasible to CP Model~4. Nevertheless, IBM ILOG CP Optimizer solver is able to heuristically transform this infeasible and incomplete initial guess into a feasible solution that potentially helps to prune the search space in the resolution process of CP Model~4. Table~\ref{tab5} shows the solutions found by IBM ILOG CP Optimizer solver when applied to the twenty medium-sized instances of CP Model~4 without using (already shown in Table~\ref{tab4}) and using the solution of the Incomplete CP Model as an initial guess. In the two-phases strategy, 2/3 of the two hours budget is allocated to the resolution of the Incomplete CP Model; while the remaining 1/3 is allocated to the resolution of CP Model~4 itself. In the table, ``\# branches~1'' corresponds to the resolution of the Incomplete CP Model; while ``\# branches~2'' corresponds to the resolution of CP Model~4. The CPU time reported for ``Incomplete CP Model + CP Model~4'' corresponds to the total CPU time for solving both models. Figures in the table show that, in instance~3, the two-phases strategy took longer than the ``single-phase strategy'' to close the gap. On the other hand, it closed the gaps of instances 9 and 17 (that were not closed by the single-phase strategy); and it reduced the gap, by improving the incumbent solution, in the other six instances in which the single-phase strategy was unable to prove optimality (namely, instances~6, 8, 11, 12, 14, and~20). Moreover, as depicted in Figure~\ref{fig:2cpmodels}, the two-phases strategy was able to prove optimality faster than the single-phase strategy on 11 out of the 12 instances in which both strategies proved optimality. As a whole, it can be inferred that the two-phases strategy for solving CP Model~4 is the most efficient way of solving the CP Optimizer formulation of the OPS scheduling problem. The figures reported in Table~\ref{tab2}, where the performances of IBM ILOG CPLEX and IBM ILOG CP Optimizer solvers are compared when applied to the MILP and the CP optimizer formulations of the OPS scheduling problem, respectively, correspond to the two-phases strategy applied to CP Model~4. As an illustration, Figure~\ref{fig:mops1} shows the solution to instance~1 of the set of medium-sized instances found by the two-phases strategy applied to CP Model~4.

\subsubsection{Experiments with large-sized instances} 

In this section, a set of~50 large-sized instances is considered. Instances are of the size of the real instances of the OPS scheduling problem. The $k$-th instance was generated with the following parameters: $n = 11 + \lceil \frac{k}{100} \times 189 \rceil$, $o_{\min} = 5$, $o_{\max} = 6 + \lceil \frac{k}{100} \times 14 \rceil$, $m_{\min} = 9 + \lceil \frac{k}{100} \times 20 \rceil$, $m_{\max} = 10 + \lceil \frac{k}{100} \times 90 \rceil$, and $q = 8$. Table~\ref{tab6} shows the main characteristics of each instance. Note that the last instance in the set has 55 machines, 250 unavailability periods, 106 jobs, 978 operations, and 1{,}581 precedence constraints; while its CP formulation has 85{,}463 integer variables and 251{,}649 constraints. Table~\ref{tab7} shows the solutions found by and the performance of IBM ILOG CP Optimizer applied to CP Model~4 and to ``Incomplete CP Model + CP Model 4''. In the table, only the number of branches is shown, since in all case the CPU time limit of two hours was reached. (Recall that, in the case in which the two phases strategy is considered, 2/3 of the time is devoted to solve the Incomplete CP Model and the remaining 1/3 of the time is used to solve CP Model~4 using as starting guess the solution to the Incomplete CP Model.) Figures in the table show that, as in the case of the medium-sized instances, the two phases strategy is more effective, in the sense that it is able to found better quality solutions in 41 out of the 50 considered instances. Figures in the table also show that in all cases a feasible solution is found and that the average gap is around~33\%. Table~\ref{tab8} shows the solution found by IBM ILOG CP Optimizer applied to ``Incomplete CP Model + CP Model 4'' with a CPU time limit of 5 minutes, 30 minutes, 2 hours, and 10 hours. The table shows that the average gap for the different CPU time limits is 41\%, 34\%, 33\%, and 31\%, respectively. These figures suggest that ``good quality'' solutions are found relatively quickly, that most of the time is used to close the gap, and that, when trying to close the gap, slightly better solutions might be found. This means that IBM ILOG CP Optimizer could be used in practice to find feasible solutions to real instances of the OPS scheduling problem. On the other hand, the non-null gap leaves space for the development of ad-hoc heuristic methods that, hopefully, would be able to find better quality solutions.

\begin{table}[ht!]
\centering 
\caption{Main features of the considered fifty large-sized instances.}
\label{tab6}
\resizebox{0.6\textwidth}{!}{\begin{tabular}{ccccccccc}
\toprule
& \multicolumn{6}{c}{Main instance characteristics} & \multicolumn{2}{c}{CP Optimizer formulation} \\
\midrule
\multirow{2}{*}{Instance} & 
\multirow{2}{*}{$m$} & 
\multirow{2}{*}{$\sum_{k=1}^m q_k$} & 
\multirow{2}{*}{$n$} & 
\multirow{2}{*}{$o$} & 
\multirow{2}{*}{$|A|$} & 
\multirow{2}{*}{$|T|$} &
\# integer & 
\multirow{2}{*}{\# constraints} \\
  &    &    &    &    &    &    & variables &   \\
\midrule  
1 & 10 & 56 & 13 & 79 & 95 & 0 & 1293 & 3792 \\ \rowcolor{gray!10.} 
2 & 10 & 44 & 15 & 98 & 120 & 0 & 1595 & 4648 \\ 
3 & 13 & 47 & 17 & 107 & 131 & 0 & 2359 & 6888 \\ \rowcolor{gray!10} 
4 & 11 & 69 & 19 & 120 & 150 & 0 & 2451 & 7209 \\ 
5 & 10 & 52 & 21 & 128 & 165 & 0 & 2235 & 6540 \\ \rowcolor{gray!10} 
6 & 15 & 64 & 23 & 141 & 180 & 0 & 3820 & 11235 \\ 
7 & 14 & 71 & 25 & 170 & 229 & 0 & 4120 & 12054 \\ \rowcolor{gray!10} 
8 & 13 & 58 & 27 & 193 & 257 & 0 & 4655 & 13631 \\ 
9 & 16 & 68 & 29 & 207 & 308 & 1 & 5592 & 16522 \\ \rowcolor{gray!10} 
10 & 19 & 81 & 30 & 201 & 266 & 0 & 6537 & 19152 \\ 
11 & 14 & 68 & 32 & 240 & 337 & 1 & 5938 & 17404 \\ \rowcolor{gray!10} 
12 & 15 & 67 & 34 & 232 & 310 & 0 & 6124 & 18115 \\ 
13 & 18 & 82 & 36 & 257 & 350 & 0 & 7860 & 23081 \\ \rowcolor{gray!10} 
14 & 13 & 60 & 38 & 259 & 350 & 0 & 5883 & 17357 \\ 
15 & 15 & 53 & 40 & 298 & 430 & 1 & 7708 & 22820 \\ \rowcolor{gray!10} 
16 & 17 & 88 & 42 & 311 & 429 & 0 & 9140 & 26801 \\ 
17 & 22 & 98 & 44 & 341 & 505 & 1 & 12788 & 37473 \\ \rowcolor{gray!10} 
18 & 13 & 57 & 46 & 345 & 486 & 1 & 7853 & 23212 \\ 
19 & 24 & 108 & 47 & 362 & 524 & 0 & 14367 & 42184 \\ \rowcolor{gray!10} 
20 & 28 & 136 & 49 & 383 & 555 & 1 & 17662 & 52005 \\ 
21 & 14 & 69 & 51 & 386 & 557 & 0 & 9544 & 28138 \\ \rowcolor{gray!10} 
22 & 20 & 79 & 53 & 398 & 572 & 0 & 12750 & 37494 \\ 
23 & 28 & 129 & 55 & 421 & 638 & 0 & 19538 & 57451 \\ \rowcolor{gray!10} 
24 & 15 & 61 & 57 & 437 & 592 & 1 & 11299 & 33312 \\ 
25 & 29 & 145 & 59 & 507 & 810 & 1 & 24190 & 70876 \\ \rowcolor{gray!10} 
26 & 26 & 120 & 61 & 487 & 743 & 0 & 21193 & 62377 \\ 
27 & 21 & 93 & 63 & 498 & 746 & 0 & 17571 & 51678 \\ \rowcolor{gray!10} 
28 & 20 & 66 & 64 & 533 & 784 & 2 & 17039 & 50072 \\ 
29 & 19 & 59 & 66 & 547 & 833 & 0 & 17555 & 51928 \\ \rowcolor{gray!10} 
30 & 26 & 130 & 68 & 563 & 831 & 0 & 24305 & 71287 \\ 
31 & 37 & 148 & 70 & 601 & 967 & 1 & 35633 & 105024 \\ \rowcolor{gray!10} 
32 & 39 & 154 & 72 & 593 & 935 & 0 & 37521 & 110493 \\ 
33 & 20 & 90 & 74 & 617 & 931 & 0 & 20027 & 58960 \\ \rowcolor{gray!10} 
34 & 35 & 161 & 76 & 656 & 1011 & 0 & 37094 & 108969 \\ 
35 & 17 & 71 & 78 & 675 & 1049 & 1 & 19362 & 57160 \\ \rowcolor{gray!10} 
36 & 19 & 81 & 80 & 700 & 1127 & 1 & 22642 & 66713 \\ 
37 & 43 & 197 & 81 & 789 & 1290 & 1 & 55235 & 162562 \\ \rowcolor{gray!10} 
38 & 43 & 186 & 83 & 706 & 1143 & 1 & 49300 & 145476 \\ 
39 & 32 & 137 & 85 & 765 & 1242 & 1 & 40211 & 118131 \\ \rowcolor{gray!10} 
40 & 26 & 93 & 87 & 812 & 1305 & 0 & 34800 & 102295 \\ 
41 & 27 & 135 & 89 & 834 & 1383 & 0 & 36985 & 109069 \\ \rowcolor{gray!10} 
42 & 38 & 152 & 91 & 821 & 1282 & 2 & 50074 & 147194 \\ 
43 & 26 & 93 & 93 & 904 & 1424 & 1 & 38750 & 114131 \\ \rowcolor{gray!10} 
44 & 18 & 89 & 95 & 921 & 1564 & 0 & 28031 & 82872 \\ 
45 & 34 & 164 & 97 & 965 & 1628 & 0 & 53196 & 156446 \\ \rowcolor{gray!10} 
46 & 28 & 127 & 98 & 868 & 1367 & 1 & 39853 & 117804 \\ 
47 & 19 & 87 & 100 & 961 & 1609 & 1 & 30667 & 90453 \\ \rowcolor{gray!10} 
48 & 53 & 228 & 102 & 942 & 1540 & 0 & 78706 & 232310 \\ 
49 & 23 & 101 & 104 & 959 & 1513 & 0 & 36885 & 108542 \\ \rowcolor{gray!10} 
50 & 55 & 250 & 106 & 978 & 1581 & 0 & 85463 & 251649 \\
\bottomrule
\end{tabular}}
\end{table}

\begin{table}[ht!]
\centering 
\caption{Description of the solutions found and effort measurements of the IBM ILOG CP Optimizer applied to the fifty large-sized instances. The left half of the table corresponds to CP Model~4 (already shown in Table~\ref{tab4}); while the right half of the table corresponds to solving first the Incomplete CP Model and passing its solution as an initial guess for the resolution of CP Model~4.}
\label{tab7}
\resizebox{0.75\textwidth}{!}{\begin{tabular}{cccccccc}
\toprule
\multirow{2}{*}{} & 
\multicolumn{3}{c}{CP Model 4} & 
\multicolumn{4}{c}{Incomplete CP Model + CP Model 4} \\
\cmidrule{2-8}
& Makespan & gap (\%) & \# branches & makespan & gap (\%) & \# branches 1 & \# branches 2 \\ 
\midrule 
1 & [387, 528] & 26.7 & 49798163 & [387, 527] & 26.6 & 120241891 & 14581110 \\ \rowcolor{gray!10}
2 & [494, 661] & 25.3 & 39711504 & [494, 650] & 24.0 & 110885682 & 13313040 \\
3 & [452, 631] & 28.4 & 39468386 & [452, 633] & 28.6 & 97607050 & 11604473 \\ \rowcolor{gray!10}
4 & [562, 763] & 26.3 & 35944769 & [562, 756] & 25.7 & 93096578 & 9904719 \\
5 & [625, 855] & 26.9 & 39251314 & [625, 828] & 24.5 & 97073106 & 10349840 \\ \rowcolor{gray!10}
6 & [490, 717] & 31.7 & 32448937 & [490, 715] & 31.5 & 85024768 & 9228219 \\
7 & [659, 943] & 30.1 & 28100210 & [659, 923] & 28.6 & 72857628 & 7788597 \\ \rowcolor{gray!10}
8 & [739, 1055] & 30.0 & 23528650 & [739, 1039] & 28.9 & 67355598 & 6681507 \\
9 & [654, 983] & 33.5 & 24454758 & [654, 980] & 33.3 & 65804340 & 6537839 \\ \rowcolor{gray!10}
10 & [547, 827] & 33.9 & 25762983 & [547, 790] & 30.8 & 63801015 & 6740620 \\
11 & [857, 1244] & 31.1 & 19205336 & [857, 1230] & 30.3 & 57729625 & 7599655 \\ \rowcolor{gray!10}
12 & [838, 1225] & 31.6 & 21224292 & [838, 1180] & 29.0 & 61503923 & 6541171 \\
13 & [699, 1073] & 34.9 & 20602257 & [699, 1011] & 30.9 & 54006978 & 5190138 \\ \rowcolor{gray!10}
14 & [1052, 1539] & 31.6 & 18026217 & [1052, 1486] & 29.2 & 50939574 & 4958123 \\
15 & [975, 1478] & 34.0 & 16007825 & [975, 1445] & 32.5 & 42382885 & 4677113 \\ \rowcolor{gray!10}
16 & [924, 1488] & 37.9 & 13721997 & [924, 1382] & 33.1 & 38344783 & 3879028 \\
17 & [763, 1150] & 33.7 & 11518247 & [763, 1117] & 31.7 & 32621308 & 2787395 \\ \rowcolor{gray!10}
18 & [1415, 2040] & 30.6 & 11688840 & [1415, 1897] & 25.4 & 36621269 & 3282986 \\
19 & [737, 1107] & 33.4 & 12841773 & [737, 1028] & 28.3 & 31061348 & 3454907 \\ \rowcolor{gray!10}
20 & [671, 1093] & 38.6 & 7883338 & [671, 1053] & 36.3 & 26628577 & 3368228 \\
21 & [1378, 2043] & 32.6 & 9278464 & [1378, 1956] & 29.6 & 38050942 & 3605124 \\ \rowcolor{gray!10}
22 & [985, 1563] & 37.0 & 11218731 & [985, 1496] & 34.2 & 37444198 & 2991893 \\
23 & [762, 1185] & 35.7 & 8102079 & [762, 1146] & 33.5 & 29827277 & 3268703 \\ \rowcolor{gray!10}
24 & [1377, 2162] & 36.3 & 7139492 & [1377, 2010] & 31.5 & 37578271 & 2841939 \\
25 & [892, 1397] & 36.1 & 4160924 & [892, 1367] & 34.7 & 20504018 & 1181601 \\ \rowcolor{gray!10}
26 & [880, 1463] & 39.8 & 7210391 & [880, 1362] & 35.4 & 26867887 & 1847124 \\
27 & [1246, 1856] & 32.9 & 8158263 & [1246, 1790] & 30.4 & 29586547 & 2630271 \\ \rowcolor{gray!10}
28 & [1396, 2175] & 35.8 & 7900895 & [1396, 2089] & 33.2 & 27240976 & 1642209 \\
29 & [1452, 2241] & 35.2 & 4952743 & [1452, 2199] & 34.0 & 28329290 & 1889844 \\ \rowcolor{gray!10}
30 & [1116, 1692] & 34.0 & 8890534 & [1116, 1769] & 36.9 & 21953739 & 1703439 \\
31 & [822, 1284] & 36.0 & 4139996 & [822, 1218] & 32.5 & 21050620 & 1652157 \\ \rowcolor{gray!10}
32 & [776, 1204] & 35.5 & 3423918 & [776, 1151] & 32.6 & 19731639 & 2487825 \\
33 & [1469, 2414] & 39.1 & 4898465 & [1469, 2276] & 35.5 & 24557450 & 2204832 \\ \rowcolor{gray!10}
34 & [951, 1490] & 36.2 & 5776499 & [951, 1483] & 35.9 & 16280020 & 1379212 \\
35 & [1932, 3161] & 38.9 & 5417895 & [1932, 3049] & 36.6 & 25096957 & 1948660 \\ \rowcolor{gray!10}
36 & [1788, 2876] & 37.8 & 4898721 & [1788, 2826] & 36.7 & 23458076 & 379399 \\
37 & [924, 1426] & 35.2 & 2985158 & [924, 1468] & 37.1 & 11587937 & 821564 \\ \rowcolor{gray!10}
38 & [832, 1281] & 35.1 & 5009442 & [832, 1192] & 30.2 & 14892180 & 1195237 \\
39 & [1214, 1895] & 35.9 & 2483387 & [1214, 1845] & 34.2 & 16606096 & 1193648 \\ \rowcolor{gray!10}
40 & [1552, 2439] & 36.4 & 3344183 & [1552, 2439] & 36.4 & 16202692 & 1411237 \\
41 & [1587, 2553] & 37.8 & 3239079 & [1587, 2524] & 37.1 & 15517269 & 873936 \\ \rowcolor{gray!10}
42 & [1111, 1905] & 41.7 & 1837052 & [1111, 1715] & 35.2 & 14722059 & 1154822 \\
43 & [1737, 2779] & 37.5 & 4525245 & [1737, 2767] & 37.2 & 16091655 & 1129221 \\ \rowcolor{gray!10}
44 & [2587, 4007] & 35.4 & 3260320 & [2587, 3908] & 33.8 & 13253041 & 1121671 \\
45 & [1446, 2331] & 38.0 & 5837368 & [1446, 2301] & 37.2 & 14297375 & 1145106 \\ \rowcolor{gray!10}
46 & [1539, 2432] & 36.7 & 3434877 & [1539, 2446] & 37.1 & 15226438 & 1126198 \\
47 & [2518, 3928] & 35.9 & 2421741 & [2518, 4040] & 37.7 & 14957027 & 961575 \\ \rowcolor{gray!10}
48 & [896, 1368] & 34.5 & 925320 & [896, 1473] & 39.2 & 6202106 & 895555 \\
49 & [2108, 3187] & 33.9 & 3352348 & [2108, 3233] & 34.8 & 17104498 & 979735 \\ \rowcolor{gray!10}
50 & [931, 1487] & 37.4 & 1857144 & [931, 1505] & 38.1 & 13503241 & 927131 \\
\midrule
Mean & 1692.48 & 34.4 & 12825409 & 1654.26 & 32.8 & 38666188 & 3821191 \\ 
\bottomrule
\end{tabular}}
\end{table}  

\begin{table}[ht!]
\centering 
\caption{Description of the solutions found by the IBM ILOG CP Optimizer applied to the fifty large-sized instances, considering the Incomplete CP Model + CP Model~4 strategy, and with increasing CPU time limits.}
\label{tab8}
\resizebox{0.75\textwidth}{!}{
\begin{tabular}{ccccccccc}
\toprule
\multirow{2}{*}{} & 
\multicolumn{2}{c}{5 minutes} & \multicolumn{2}{c}{30 minutes} & 
\multicolumn{2}{c}{2 hours} & \multicolumn{2}{c}{10 hours} \\
\cmidrule{2-9}
& Makespan & gap (\%) & Makespan & gap (\%) & Makespan & gap (\%) & Makespan & gap (\%) \\
\midrule  
1 & [387, 538] & 28.1 & [387, 530] & 27.0 & [387, 527] & 26.6 & [387, 521] & 25.7\\ \rowcolor{gray!10}
2 & [494, 663] & 25.5 & [494, 654] & 24.5 & [494, 650] & 24.0 & [494, 649] & 23.9\\
3 & [452, 653] & 30.8 & [452, 635] & 28.8 & [452, 633] & 28.6 & [452, 623] & 27.4\\ \rowcolor{gray!10}
4 & [562, 780] & 27.9 & [562, 755] & 25.6 & [562, 756] & 25.7 & [562, 756] & 25.7\\
5 & [625, 860] & 27.3 & [625, 837] & 25.3 & [625, 828] & 24.5 & [625, 825] & 24.2\\ \rowcolor{gray!10}
6 & [490, 724] & 32.3 & [490, 718] & 31.8 & [490, 715] & 31.5 & [490, 698] & 29.8\\
7 & [659, 964] & 31.6 & [659, 938] & 29.7 & [659, 923] & 28.6 & [659, 902] & 26.9\\ \rowcolor{gray!10}
8 & [739, 1091] & 32.3 & [739, 1044] & 29.2 & [739, 1039] & 28.9 & [739, 1022] & 27.7\\
9 & [654, 1019] & 35.8 & [654, 966] & 32.3 & [654, 980] & 33.3 & [654, 936] & 30.1\\ \rowcolor{gray!10}
10 & [547, 902] & 39.4 & [547, 900] & 39.2 & [547, 790] & 30.8 & [547, 781] & 30.0\\
11 & [857, 1290] & 33.6 & [857, 1231] & 30.4 & [857, 1230] & 30.3 & [857, 1182] & 27.5\\ \rowcolor{gray!10}
12 & [838, 1257] & 33.3 & [838, 1223] & 31.5 & [838, 1180] & 29.0 & [838, 1178] & 28.9\\
13 & [699, 1084] & 35.5 & [699, 1010] & 30.8 & [699, 1011] & 30.9 & [699, 998] & 30.0\\ \rowcolor{gray!10}
14 & [1052, 1557] & 32.4 & [1052, 1514] & 30.5 & [1052, 1486] & 29.2 & [1052, 1449] & 27.4\\
15 & [975, 1542] & 36.8 & [975, 1476] & 33.9 & [975, 1445] & 32.5 & [975, 1427] & 31.7\\ \rowcolor{gray!10}
16 & [924, 1516] & 39.1 & [924, 1420] & 34.9 & [924, 1382] & 33.1 & [924, 1343] & 31.2\\
17 & [763, 1131] & 32.5 & [763, 1080] & 29.4 & [763, 1117] & 31.7 & [763, 1037] & 26.4\\ \rowcolor{gray!10}
18 & [1415, 2014] & 29.7 & [1415, 1918] & 26.2 & [1415, 1897] & 25.4 & [1415, 1898] & 25.4\\
19 & [737, 1236] & 40.4 & [737, 1046] & 29.5 & [737, 1028] & 28.3 & [737, 1008] & 26.9\\ \rowcolor{gray!10}
20 & [671, 1135] & 40.9 & [671, 1132] & 40.7 & [671, 1053] & 36.3 & [671, 1050] & 36.1\\
21 & [1378, 2104] & 34.5 & [1378, 1992] & 30.8 & [1378, 1956] & 29.6 & [1378, 1939] & 28.9\\ \rowcolor{gray!10}
22 & [985, 1639] & 39.9 & [985, 1642] & 40.0 & [985, 1496] & 34.2 & [985, 1573] & 37.4\\
23 & [762, 1336] & 43.0 & [762, 1128] & 32.4 & [762, 1146] & 33.5 & [762, 1082] & 29.6\\ \rowcolor{gray!10}
24 & [1377, 2135] & 35.5 & [1377, 2073] & 33.6 & [1377, 2010] & 31.5 & [1377, 1942] & 29.1\\
25 & [892, 1524] & 41.5 & [892, 1442] & 38.1 & [892, 1367] & 34.7 & [892, 1241] & 28.1\\ \rowcolor{gray!10}
26 & [880, 1581] & 44.3 & [880, 1486] & 40.8 & [880, 1362] & 35.4 & [880, 1348] & 34.7\\
27 & [1246, 1833] & 32.0 & [1246, 1791] & 30.4 & [1246, 1790] & 30.4 & [1246, 1741] & 28.4\\ \rowcolor{gray!10}
28 & [1396, 2185] & 36.1 & [1396, 2171] & 35.7 & [1396, 2089] & 33.2 & [1396, 2003] & 30.3\\
29 & [1452, 2256] & 35.6 & [1452, 2298] & 36.8 & [1452, 2199] & 34.0 & [1452, 2282] & 36.4\\ \rowcolor{gray!10}
30 & [1116, 2473] & 54.9 & [1116, 1590] & 29.8 & [1116, 1769] & 36.9 & [1116, 1561] & 28.5\\
31 & [822, 1296] & 36.6 & [822, 1295] & 36.5 & [822, 1218] & 32.5 & [822, 1288] & 36.2\\ \rowcolor{gray!10}
32 & [434, 1214] & 64.3 & [776, 1208] & 35.8 & [776, 1151] & 32.6 & [776, 1165] & 33.4\\
33 & [1469, 2698] & 45.6 & [1469, 2398] & 38.7 & [1469, 2276] & 35.5 & [1469, 2166] & 32.2\\ \rowcolor{gray!10}
34 & [409, 1545] & 73.5 & [951, 1604] & 40.7 & [951, 1483] & 35.9 & [951, 1403] & 32.2\\
35 & [1932, 3145] & 38.6 & [1932, 3099] & 37.7 & [1932, 3049] & 36.6 & [1932, 2892] & 33.2\\ \rowcolor{gray!10}
36 & [1788, 3167] & 43.5 & [1788, 2819] & 36.6 & [1788, 2826] & 36.7 & [1788, 2691] & 33.6\\
37 & [924, 1463] & 36.8 & [924, 1671] & 44.7 & [924, 1468] & 37.1 & [924, 1350] & 31.6\\ \rowcolor{gray!10}
38 & [832, 1264] & 34.2 & [832, 1200] & 30.7 & [832, 1192] & 30.2 & [832, 1144] & 27.3\\
39 & [452, 1880] & 76.0 & [1214, 1870] & 35.1 & [1214, 1845] & 34.2 & [1214, 1786] & 32.0\\ \rowcolor{gray!10}
40 & [1552, 2463] & 37.0 & [1552, 2509] & 38.1 & [1552, 2439] & 36.4 & [1552, 2347] & 33.9\\
41 & [1587, 2550] & 37.8 & [1587, 2522] & 37.1 & [1587, 2524] & 37.1 & [1587, 2542] & 37.6\\ \rowcolor{gray!10}
42 & [618, 1732] & 64.3 & [1111, 1688] & 34.2 & [1111, 1715] & 35.2 & [1111, 1598] & 30.5\\
43 & [529, 2784] & 81.0 & [1737, 2779] & 37.5 & [1737, 2767] & 37.2 & [1737, 2599] & 33.2\\ \rowcolor{gray!10}
44 & [2587, 4030] & 35.8 & [2587, 4063] & 36.3 & [2587, 3908] & 33.8 & [2587, 3770] & 31.4\\
45 & [1446, 2378] & 39.2 & [1446, 2466] & 41.4 & [1446, 2301] & 37.2 & [1446, 2210] & 34.6\\ \rowcolor{gray!10}
46 & [1539, 2439] & 36.9 & [1539, 2497] & 38.4 & [1539, 2446] & 37.1 & [1539, 2293] & 32.9\\
47 & [2518, 4195] & 40.0 & [2518, 3757] & 33.0 & [2518, 4040] & 37.7 & [2518, 3824] & 34.2\\ \rowcolor{gray!10}
48 & [544, 1377] & 60.5 & [896, 1506] & 40.5 & [896, 1473] & 39.2 & [896, 1290] & 30.5\\
49 & [2108, 4065] & 48.1 & [2108, 3012] & 30.0 & [2108, 3233] & 34.8 & [2108, 2897] & 27.2\\ \rowcolor{gray!10}
50 & [931, 2244] & 58.5 & [931, 1570] & 40.7 & [931, 1505] & 38.1 & [931, 1478] & 37.0 \\ 
\midrule
Mean & 1779.02 & 41.0 & 1683.46 & 34.1 & 1654.26 & 32.8 & 1594.56 & 30.6 \\
\bottomrule
\end{tabular}}
\end{table}  

\section{Conclusions and future works}\label{sec:concl}

In this work, a challenging real scheduling problem with several complicating features was introduced. The problem, named OPS scheduling problem, can be seen as a flexible job shop scheduling problem with sequence flexibility. MILP and CP Optimizer formulations of the problem were presented and analyzed. The capacity of IBM ILOG CPLEX MILP and IBM ILOG CP Optimizer solvers for dealing with the proposed formulations was investigated. While the MILP solver is a general-purpose solver; the CP Optimizer solver was born as a solver dedicated to scheduling problems, with its own modeling language that fully explores the structure of the underlying problem. Thus, it is not a surprise the latter to outperform the former by a large extent in the numerical experiments with the OPS scheduling problem presented in this work. The obtained result is in agreement with a previous comparison of similar nature presented in~\cite{vilim2015failure} and~\cite{Laborie2018a} for the job shop and the flexible job shop scheduling problems. 

The problem under consideration corresponds to a real problem of a European company that runs several online printing shops all over the Continent. These shops process, on average, $20{,}000$ orders per day and they recompute the scheduling considering recently arrived orders several times per day. The OPS scheduling problem is an NP-hard problem and the CP Optimizer mixes exact and heuristic strategies that were not specifically devised for the problem at hand. However, numerical experiments in this work show that the CP Optimizer is able to find feasible solution to large-sized instances; thus being an alternative to tackle the OPS scheduling problem in practice. On the other hand, it is not expected the CP Optimizer to be competitive with ad-hoc heuristics fully exploiting the specificities of the problem. These facts suggest that the development of heuristic methods to deal with the considered problem is a promising alternative; and this will be the subject of future work. 

\section*{Acknowledgement}
\noindent
The experiments presented in this paper were carried out using the HPC facilities of the University of Luxembourg~\cite{VBCG_HPCS14} {\small -- see \url{https://hpc.uni.lu}}. This work has been partially supported by the Brazilian funding agencies FAPESP (grants 2013/07375-0, 2016/01860-1, and 2018/24293-0) and CNPq (grants 306083/2016-7 and 302682/2019-8). The authors would like to thank the careful reading and the comments of the reviewers that helped a lot to improve the quality of this work.

\section*{References}

\bibliography{ops1}

\end{document}